  \documentclass[11pt,a4paper,twoside]{article}

  \usepackage[margin=3.07cm]{geometry}
  \usepackage{amsmath,amssymb,amsthm,amsfonts,bbding,stmaryrd,dsfont,wasysym,pifont,xspace,xcolor}
  \usepackage{parskip,fancyhdr,url}
  \usepackage{fourier}
  \usepackage{enumitem}
  \usepackage{mathrsfs}
  \usepackage{todonotes}
  \usepackage[
     breaklinks,colorlinks,
     citecolor=blue,linkcolor=blue,urlcolor=teal
  ]{hyperref}
  \usepackage{tikz}
  \usetikzlibrary{
    arrows,
    calc
    }

  \begingroup
    \makeatletter
    \@for\theoremstyle:=definition,remark,plain\do{%
      \expandafter\g@addto@macro\csname th@\theoremstyle\endcsname{%
        \addtolength\thm@preskip\parskip
        }%
      }
  \endgroup

  \makeatletter
    \def\tagform@#1{\maketag@@@{%
     \textbf{(\ignorespaces#1\unskip\@@italiccorr)}}}%
     \renewcommand{\eqref}[1]{\textup{\maketag@@@{(\ignorespaces%
          {\ref{#1}}\unskip\@@italiccorr)}}}
  \makeatother

  \newcommand\address[1]{}
  \newcommand\email[1]{}
  \newcommand\dedicatory[1]{}

  \theoremstyle{plain}
  \newtheorem{theorem}[equation]{Theorem}
  \newtheorem{proposition}[equation]{Proposition}
  \newtheorem{corollary}[equation]{Corollary}
  \newtheorem{lemma}[equation]{Lemma}

  \theoremstyle{definition}
  \newtheorem{definition}[equation]{Definition}
  \newtheorem{remark}[equation]{Remark}

  \theoremstyle{axiom}

  \numberwithin{equation}{section}
  
  \theoremstyle{problem}

  \newcommand{\set}[1]{\ensuremath{\left\{ {#1} \right\}}\xspace}

  \newcommand{\st}{\ensuremath{\,\, \colon \,\,}\xspace} 
  \newcommand{\from}{\ensuremath{\colon \thinspace}\xspace} 

  \newcommand{\halmos}{\hfill$\square$}
  
  \DeclareMathOperator{\area}{area}
  
  \DeclareMathOperator{\width}{width}

  \DeclareMathOperator{\hol}{hol}
  \DeclareMathOperator{\diam}{diam}

  \DeclareMathOperator{\Aff}{Aff}

  \DeclareMathOperator{\Aut}{Aut}

  \newcommand{\sys}{\ensuremath{\boldsymbol{\sigma}}\xspace}
  \newcommand{\smark}{\ensuremath{\boldsymbol{\mu}}\xspace}

  \newcommand{\Z}{\ensuremath{\mathbb{Z}}\xspace}
  
  \newcommand{\R}{\ensuremath{\mathbb{R}}\xspace}
  \newcommand{\CC}{\ensuremath{\mathbb{C}}\xspace}
  \newcommand{\SL}{\ensuremath{\mathrm{SL}(2,\R)}\xspace}
  \newcommand{\PSL}{\ensuremath{\mathrm{PSL}\xspace}}
  \newcommand{\RP}{\R \mathrm{P}^1}
  \newcommand{\Hy}{\ensuremath{\mathbb{H}}\xspace}
  \newcommand{\Ann}{\ensuremath{\mathbb{A}}\xspace}

  \newcommand{\T}{\ensuremath{\mathcal{T}}\xspace} 
   
  \newcommand{\mcg}{\ensuremath{\text{MCG}}\xspace}  
   
  \newcommand{\PMF}{\ensuremath{\mathcal{PMF}}\xspace} 
  \newcommand{\QD}{\ensuremath{\text{QD}}\xspace} 
   
  \newcommand{\C}{\ensuremath{\mathcal{C}}\xspace} 
   
  \renewcommand{\AC}{\ensuremath{\mathcal{AC}}\xspace} 
  \newcommand{\ParC}{\ensuremath{\mathrm{PCyl}}\xspace} 
  
   \newcommand{\Par}{\ensuremath{\mathcal{P}}\xspace} 
  
  \newcommand{\G}{\ensuremath{\mathcal{G}}\xspace} 
   
    \newcommand{\Ext}{\ensuremath{\mathrm{Ext}}\xspace} 
       
	\newcommand{\X}{\ensuremath{\mathcal{X}}\xspace} 
	\newcommand{\Y}{\ensuremath{\mathcal{Y}}\xspace} 
	
  \newcommand{\bG}{\ensuremath{\mathbf{G}}\xspace} 
       
   \newcommand{\Disc}{\ensuremath{\mathbf{H}}\xspace} 
   
   \newcommand{\Core}{\ensuremath{\mathrm{N}}\xspace}
   \newcommand{\elCore}{\ensuremath{\mathrm{N}^{el}}\xspace}
   \newcommand{\trCore}{\ensuremath{\mathrm{N}^{tr}}\xspace}
   \newcommand{\dist}{\ensuremath{\mathrm{d}}\xspace}
   \newcommand{\eldist}{\ensuremath{\mathrm{d}^{el}}\xspace}
   \newcommand{\trdist}{\ensuremath{\mathrm{d}^{tr}}\xspace}
        
   \newcommand{\elTeich}{\ensuremath{\T^{el}}\xspace}

   \newcommand{\tw}{\ensuremath{\mathrm{tw}}\xspace}
   
   \newcommand{\Mark}{\ensuremath{\mathcal{M}}\xspace} 

  \newcommand{\sA}{\mathsf{A}}
  \newcommand{\sB}{\mathsf{B}}
  \newcommand{\sK}{\mathsf{K}}

  \newcommand{\sN}{\mathsf{N}}
  
  \newcommand{\sY}{\mathsf{Y}}

  \newcommand{\sr}{\mathsf{r}}

  \newcommand{\proofof}[1]{\hfill\newline\noindent\emph{Proof of {#1}.} }
  
  \newcommand{\param}{{\mathchoice{\mkern1mu\mbox{\raise2.2pt\hbox{$
  \centerdot$}}
  \mkern1mu}{\mkern1mu\mbox{\raise2.2pt\hbox{$\centerdot$}}\mkern1mu}{
  \mkern1.5mu\centerdot\mkern1.5mu}{\mkern1.5mu\centerdot\mkern1.5mu}}}

  \fancypagestyle{plain}
  
  \AtEndDocument{\bigskip{\footnotesize%
  Topology and Geometry of Manifolds Unit, Okinawa Institute of Science and Technology Graduate University,
  1919-1 Tancha, Onna-son, Kunigami-gun, Okinawa, 904-0495, Japan
  \par
  \texttt{robert.tang@oist.jp} 
  }}

\begin{document}

  \title    {Affine diffeomorphism groups are undistorted}
    \author   {Robert Tang} \date{13 December 2019}

  \maketitle \thispagestyle{empty}

  \begin{abstract} 

    The affine diffeomorphism group $\Aff(S,q)$ of a half-translation surface $(S,q)$
    comprise the self-diffeomorphisms with constant differential away from the singularities.
    This group coincides with the stabiliser of the associated Teichm\"uller disc
    under the action of the mapping class group on Teichm\"uller space.
    We prove that any finitely generated subgroup of $\Aff(S,q)$ is undistorted in the mapping class group.
    We also show that the systole map restricted to the associated electrified Nielsen core in the
    Teichm\"uller disc is a quasi-isometric embedding into the curve complex.
    
  \end{abstract}

\section{Introduction}

 For a finite type surface $S$, the mapping class group $\mcg(S)$ is the (finitely generated) group
 of orientation-preserving self-homeomorphisms of $S$ up to isotopy.
 There has been considerable interest in understanding the large scale geometry of mapping class groups.
 In particular, determining the distortion properties of naturally occurring subgroups plays a central role in this regard;
 see Problem 3.7 in \cite{Farb-problems}.
 On the one hand, subsurface mapping class groups \cite{MM2, Hamen-mcg1}
 and convex cocompact subgroups \cite{farb-convcoc} are undistorted, while on the other, handlebody groups \cite{HH-handle} and Torelli groups \cite{BFP-distortion} have been shown to be distorted.
 
 Our focus is on the distortion of subgroups that stabilise Teichm\"uller discs.
 Recall that the Teichm\"uller space $\T(S)$ parameterises the marked hyperbolic structures on $S$ up to isotopy.
 The mapping class group acts by isometries on $\T(S)$ with respect to the Teichm\"uller metric.
 A Teichm\"uller disc is a geodesically embedded copy of the hyperbolic plane in $\T(S)$
 that arises from the $\SL$--orbit of a quadratic differential; 
 its stabiliser under the action of $\mcg(S)$ on $\T(S)$
 can be naturally identified with the \emph{affine diffeomorphism group} of any quadratic differential generating the given Teichm\"uller disc.  These form an important class of groups, and are closely related to the study of billiard dynamics and translation surfaces \cite{Veech-dich, MT-billiards, GHS-affine}.
 
 \begin{theorem}\label{thm:main}
  Any finitely generated subgroup of $\mcg(S)$ stabilising a Teichm\"uller disc in $\T(S)$ is undistorted.
  In particular, affine diffeomorphism groups of Veech surfaces are undistorted.
 \end{theorem}
 
 \begin{remark}
  In the cases where $\mcg(S) \cong \mathrm{SL}(2,\Z)$, this is a consequence of the fact that finitely generated subgroups of a virtually free group are quasiconvex, and hence undistorted.
 \end{remark}
 
 \begin{remark}
  The subgroups under consideration in Thereom \ref{thm:main} are virtually free. However, there exist distorted free subgroups of $\mcg(S)$; see for example, the \emph{point-pushing subgroups} in the case where $S$ is a punctured surface \cite{BFP-distortion}.
 \end{remark}

 Much of this paper is devoted to establishing bounded geometry results for Nielsen cores in Teichm\"uller discs.
 We assume throughout that $\Gamma \leq \mcg(S)$ is finitely generated subgroup stabilising a Teichm\"uller disc $\Disc(\Gamma)$. We also assume that $\Gamma$ is not virtually cyclic.
 Viewing $\Disc(\Gamma)$ as a copy of the hyperbolic plane,
 $\Gamma$ acts as a finitely generated Fuchsian group (upon passing to a finite-index quotient) \cite{Veech-dich};
 the \emph{Nielsen core} $\Core(\Gamma) \subseteq \Disc(\Gamma)$ is the convex hull of its limit set in $\partial \Disc(\Gamma)$. The finite generation assumption implies that the quotient $\Core(\Gamma)/\Gamma$ is a finite-area hyperbolic orbifold (possibly with geodesic boundary) and hence has finitely many cusps. Choose a sufficiently small horocyclic neighbourhood of each cusp so that their preimages in $\Core(\Gamma)$ give a collection of pairwise disjoint horodiscs.
 We construct two variants of the Nielsen core using this family of horodiscs.
 The \emph{electrified Nielsen core} $\elCore(\Gamma)$ is obtained from $\Core(\Gamma)$
 by forcing each horodisc to have uniformly bounded diameter  (see Section \ref{sec:qi} for the definition of an electrified space).
 The \emph{trunctated Nielsen core} $\trCore(\Gamma)$ is the complement of the interiors of all
 horodiscs in $\Core(\Gamma)$ equipped with the path metric.
 The action of $\Gamma$ on $\trCore(\Gamma)$ is geometric (properly discontinuous and cocompact),
 and so any orbit map $\Gamma \to \trCore(\Gamma)$ is a quasi-isometry by the \u Svarc--Milnor Lemma.
 
 We prove that the electrified and truncated Nielsen cores respectively quasi-isometrically embed into
 to the curve graph $\C(S)$ and marking graph $\Mark(S)$ (see Section \ref{sec:comb} for background on combinatorial complexes).
 The systole map $\sys \from \T(S) \to \C(S)$ sends a hyperbolic surface to its set of shortest curves.
 A celebrated theorem of Masur and Minsky is that the systole map is a quasi-isometry, where $\T(S)$ is equipped with
 the \emph{electrified} Teichm\"uller metric  and $\C(S)$ with the combinatorial metric \cite{MM1}.
 Leininger asks whether this still holds if the systole map is restricted to an electrified Teichm\"uller disc arising from a Veech surface (in which case the Nielsen core is the full Teichm\"uller disc). We give a positive answer in a more general setting.

 \begin{theorem}\label{thm:qi_systole}
  The restricted systole map $\sys \from \elCore(\Gamma) \to \C(S)$ is a $\Gamma$--equivariant quasi-isometric embedding.
   \end{theorem}

 As a consequence, the natural inclusion $\Core(\Gamma) \hookrightarrow \T(S)$ is a quasi-isometric
 embedding when both spaces are equipped with their respective electrified metrics (see Corollary \ref{cor:elteich}).
 
 Next, we consider the analogous statement for the truncated Nielsen core and the marking graph.
 Masur and Minsky show that the mapping class group acts geometrically on $\Mark(S)$ \cite{MM2}.
 They also define a $\mcg(S)$--equivariant \emph{short marking map} $\smark \from \T(S) \to \Mark(S)$;
 see Section \ref{sec:systole} for details.
 
 \begin{theorem}\label{thm:qi_marking}
  The restricted short marking map $\smark \from \trCore(\Gamma) \to \Mark(S)$
  is a $\Gamma$--equivariant quasi-isometric embedding.
 \end{theorem}
 
 There is an analogous consequence: the inclusion of $\trCore(\Gamma)$ into the thick part of Teichm\"uller space
 is a quasi-isometric embedding (see Corollary \ref{cor:thick}).

 \proofof{Theorem \ref{thm:main}}
  If $\Gamma$ is virtually cyclic then it is undistorted in $\mcg(S)$ \cite{FLM-rank1}, so we may assume otherwise. 
  Choose a basepoint $x_0 \in \trCore(\Gamma)$ and let $\mu_0$ be a short marking at $x_0$.
  Using the above theorem and the fact that the action of $\Gamma$ on $\trCore(\Gamma)$ is geometric, we deduce that the orbit map $\Gamma \to \Mark(S)$ given by $g \mapsto g\cdot \mu_0 = \smark(g\cdot x_0)$ is a quasi-isometric embedding.
  Since $\mcg(S)$ acts geometrically on $\Mark(S)$,
  it follows that the inclusion $\Gamma \hookrightarrow \mcg(S)$ is a quasi-isometric embedding.
 \halmos
 
 The proofs of Theorems \ref{thm:qi_systole} and \ref{thm:qi_marking} rely on the following technical result involving subsurface projections; see Section \ref{sec:comb} for the definition of the map $\pi_Y \from \C(S) \to \C(Y)$ where $Y\subseteq S$ is an essential subsurface.
 Any parabolic subgroup of $\Gamma$ has an invariant multicurve on $S$; call an annulus on $S$ \emph{parabolic} for $\Gamma$ if its core curve is a component of such an invariant multicurve.
 
 \newpage
 \begin{proposition}
  There exists a constant $D = D(S, \Gamma)$ such that for any essential subsurface ${Y \subsetneq S}$,
  the image $\pi_Y( \smark(\Core(\Gamma)))$ has infinite diameter in $\C(Y)$ if $Y$ is a parabolic annulus for $\Gamma$;
  and diameter at most $D$ in $\C(Y)$ otherwise.
 \end{proposition}
 
 This statement originally appears as Lemma 5.16 in a paper of Durham, Hagen, and Sisto \cite{HHS-boundary}, and is a key ingredient
 in proving that $\Gamma$ is hierarchically hyperbolic with respect to its parabolic subgroups.
 However, after discussions with the authors, it became apparent that there is a mistake in their proof, and so we shall give an alternative proof in Section \ref{sec:bounded}.
 (In their paper, $\Gamma$ is called a \emph{Veech subgroup} but we shall not use this term
 in order to avoid confusion with Veech groups which are subgroups of $\PSL(2,\R)$.)

 \subsection*{Acknowledgements}
 
 The author thanks Chris Leininger for interesting discussions and for asking whether Theorem \ref{thm:qi_systole} is true.
 Thanks also to Matthew Durham, Mark Hagen, and Alessandro Sisto for engaging conversations regarding their work,
 and to Richard Webb for helpful comments.
 This work was supported by a JSPS KAKENHI Grant-in-Aid for Early-Career Scientists (No.~19K14541). 

\section{Background}

 We begin by recalling some standard notions from coarse geometry.
 Given $a,b \geq 0$ and $\sK > 0$,
 write $a \prec_\sK b$ to mean $a \leq \sK\cdot a + \sK$, and
 $a \asymp_\sK b$ when $a \prec_\sK b$ and $b \prec_\sK a$.
 When the constant $\sK$ can be chosen to depend only on the topology of
 a surface $S$, we shall also write $a \prec b$ and $a \asymp b$ for simplicity.

 A map between metric spaces $f \from \X \to \Y$ is called \emph{coarsely Lipschitz}
 if there exists some $\sK > 0$ such that
 \[\dist_\Y(f(x), f(y)) \prec_\sK \dist_\X(x,y)\]
 for all $x,y\in\X$.  Furthermore, if
 \[\dist_\Y(f(x), f(y)) \asymp_\sK \dist_\X(x,y)\]
 for all $x,y\in\X$ then we call $f$ a \emph{quasi-isometric embedding}.
 A quasi-isometric embedding with coarsely dense image is called a
 \emph{quasi-isometry}.
 A \emph{quasigeodesic} is a quasi-isometric embedding of an interval.
 These notions also make sense if $f$ is multi-valued, so long as we assume that the
 image of each $x\in\X$ is non-empty and has uniformly bounded diameter.
 When $\X$ is a combinatorial complex, we shall adopt the convention
 \[\dist_\X(X,Y) := \diam_\X(X \cup Y) \]
 for sets $X,Y \subseteq \X$; this ensures that the triangle inequality holds.

 If $G$ is a finitely generated group, then a finitely generated subgroup $H \leq G$
 is \emph{undistorted} if the inclusion map $H \hookrightarrow G$ 
 is a quasi-isometric embedding with respect to (any of) their word metrics.
 
 Given $c \geq 0$, the \emph{cutoff function} is defined by
 $[t]_c = t$ if $t \geq c$, and $[t]_c = 0$ otherwise.

 \subsection{Combinatorial complexes}\label{sec:comb}

 Throughout this paper, we shall assume that $S$ is a connected, orientable surface (without boundary)
 of finite genus and with a finite set $Z$ of punctures. Furthermore, we assume that its \emph{complexity}
 \[\xi(S) := 3\cdot \textrm{genus}(S) + |Z| - 3\]
 is at least 2.
 The \emph{mapping class group} $\mcg(S)$ is the group of orientation-preserving self-homeomorphisms
 of $S$ up to isotopy.
 Mapping class groups have been fruitfully studied via their actions on various graphs associated to $S$.
 Each of these graphs is endowed with the standard combinatorial metric,
 where each edge is isometrically identified with an interval of unit length.
 For further reference, see \cite{MM1, MM2}.

 A \emph{curve} on $S$ is an (isotopy class of an) embedded loop on $S$
 that is not homotopic to a point or into a puncture.
 An \emph{arc} is (a proper isotopy class of) an interval on $S$ 
 that has embedded interior, with endpoints contained in the set of punctures,
 and cannot be homotoped into a puncture.
 The \emph{arc-and-curve graph} $\AC(S)$ has as vertices the arcs and curves on $S$,
 with edges connecting pairs of vertices whenever the corresponding arcs or curves have
 disjoint representatives.
 The \emph{curve graph} $\C(S)$ is the induced subgraph of $\AC(S)$ whose vertices
 are the curves.
 Both of these graphs are connected, locally infinite, have infinite diameter, and are Gromov hyperbolic.
 Furthermore, the inclusion map $\C(S) \hookrightarrow \AC(S)$ is a quasi-isometry.

 There is a modified definition of the arc-and-curve graph in the case of the (closed) annulus $\Ann$:
 vertices of $\AC(\Ann)$ are embedded arcs connecting the two boundary components of $\Ann$
 considered up to isotopy fixing their endpoints; while edges are defined as usual.
 The resulting graph $\AC(\Ann)$ is quasi-isometric to $\Z$. 

 Next, we consider the \emph{marking graph} $\Mark(S)$.
 We shall not recall the full definition; instead, we state some facts that suffice for our purposes.
 A \emph{marking} $\mu$ on $S$ consists of
 a pants decomposition $base(\mu)$ of $S$, called the set of \emph{base curves},
 and for each $\beta \in base(\mu)$, another curve, called a \emph{transversal}, that
 intersects $\beta$ exactly once or twice and is disjoint from all other base curves.
 Each marking on $S$ has diameter at most 3 as a subset of $\C(S)$.
 The set of markings on $S$ form the vertices of $\Mark(S)$,
 with edges defined using a rule which guarantees the following:
 \begin{itemize}
   \item $\Mark(S)$ is connected, locally finite, and admits a geometric action by $\mcg(S)$,
   \item if $\mu, \mu' \in \Mark(S)$ are adjacent, then $\diam_{\C(S)}(\mu \cup \mu') \leq 4$.
 \end{itemize}
 By the \u Svarc--Milnor Lemma, any orbit map $\mcg(S) \to \Mark(S)$ is a quasi-isometry.

 Let $Y \subseteq S$ be a non-pants essential subsurface .
 For any $\alpha\in\AC(S)$ intersecting $Y$ essentially,
 the \emph{subsurface projection} $\pi_Y(\alpha) \subseteq \AC(Y)$ is defined as follows.
 Equip $S$ with a complete hyperbolic metric (the choice of metric does not matter for this construction).
 Let $S^Y$ be the cover corresponding to $\pi_1(Y)$
 and $\tilde{\alpha}$ be the pre-image of the geodesic representative of $\alpha$.  The Gromov compactification $\overline{S^Y}$ admits a natural identification with $Y$;
 we then set $\pi_Y(\alpha)$ to be all essential arcs and curves appearing in the closure of
 $\tilde{\alpha}$ in $\overline{S^Y}$.
 We extend subsurface projections to subsets of $\AC(S)$ by taking the union
 of the images of the individual arcs and curves.
 Subsurface projections can also be defined for a (measured) foliation $F$ on $S$:
 set $\pi_Y(F)$ to be the set of all essential arcs or curves on $\overline{S^Y}$ that descend
 to a leaf of $F$. (It may be the case that $Y$ is the subsurface filled by some non-compact leaf
 of $F$, in which case $\pi_Y(F)$ is empty.)

 \begin{lemma}[\cite{MM2}]
  For any essential subsurface $Y \subseteq S$, the map   ${\pi_Y \from \Mark(S) \to \AC(Y)}$ is uniformly coarsely Lipschitz. \halmos
 \end{lemma}

 The following is the celebrated distance formula for the mapping class group.
 We shall use $\dist_Y$ as shorthand notation for $\dist_{\AC(Y)}$.
 In the case where $Y$ is an annulus with core curve $\alpha$, we also write $\dist_{\alpha}$ in place of $\dist_Y$.

 \begin{theorem}[\cite{MM2}]\label{thm:mm_dist}
  There exists a constant $c_1 = c_1(S)$ such that for all $c \geq c_1$,
  there exists $\sA_1 > 0$ such that
  \[\dist_{\Mark(S)}(\mu, \mu') \asymp_{\sA_1} \sum_Y \left[\dist_Y(\mu, \mu')\right]_c \]
  for all markings $\mu, \mu' \in \Mark(S)$,
  where the sum is taken over all essential subsurfaces $Y \subseteq S$. \halmos
 \end{theorem}

 We shall also recall a distance bounds for the arc-and-curve graph in terms of geometric intersection numbers. For closed surfaces, Hempel showed that
 \begin{equation}\label{eqn:hempel}
  \dist_S(\alpha, \beta) \leq 2 \log \iota(\alpha, \beta) + 2
 \end{equation}  
 whenever $\alpha, \beta$ are curves satisfying $\iota(\alpha, \beta) \neq 0$  \cite{Hempel-cc};
 Schleimer extended this result to all non-annular surfaces \cite{Saul-notes}.
 By a standard argument, the bound also holds (at the cost of increasing the additive constant)
 if $\alpha$ or $\beta$ are arcs. The distance between two arcs
 $\alpha$ and $\beta$ in $\AC(\Ann)$ agrees with $\iota(\alpha, \beta)$ up to a uniform additive error.

 \subsection{Teichm\"uller space, quadratic differentials, and Teichm\"uller discs}\label{sec:teich}

 We now recall some background on Teichm\"uller theory and half-translation surfaces; refer to
 \cite{FLP-travaux, Strebel-quadratic, FM-primer} for further details.
 The \emph{Teichm\"uller space} $\T(S)$ of $S$ is the space of marked complete hyperbolic metrics on
 $S$ (up to isotopy).  By the Uniformisation Theorem, this is equivalent to the space of marked conformal
 structures on $S$. The mapping class group acts on $\T(S)$ by change of marking. Moreover, $\T(S)$ is homeomorphic to $\R^{2\xi(S)}$.

 The cotangent bundle to $\T(S)$ is naturally identified with the space $\QD(S)$ of \emph{quadratic differentials} (up to isotopy).
 By a quadratic differential on $S$, we mean a Riemann surface $x \in \T(S)$ equipped with an integrable meromorphic quadratic differential that has simple poles at (and only at) the punctures.
 We shall use $q\in\QD(S)$ to denote a quadratic differential, with the underlying conformal structure implicit in the notation. Let $\QD(x)$ be the space of quadratic differentials with underlying conformal structure $x \in \T(S)$.
 A quadratic differential $q\in\QD(x)$ gives rise to \emph{natural co-ordinates}:
 these are defined in a neighbourhood of a point $z_0$ by
 \[z \mapsto \int_{z_0}^z \sqrt{q(w)} dw\]
 in terms of the complex co-ordinates from $x$.
 These natural co-ordinates give an atlas away from the zeroes of $q$ where the transition maps
 are of the form $z \mapsto \pm z + c$ for some $c \in \CC$. Pulling back the Euclidean metric
 on $\CC$ via this atlas endows $S$ with a locally Euclidean metric away from the zeroes of $q$,
 together with a preferred choice of vertical slope. The metric completion yields a singular Euclidean metric, known as a \emph{half-translation structure}, where each zero of order $p$ becomes a Euclidean cone point with cone angle $(p+2)\pi$; in particular, poles have cone angle $\pi$.
 We shall use $(S,q)$, or simply $q$, to denote $S$ equipped with this half-translation structure (with the choice of vertical slope). The integrability assumption ensures that the Euclidean area of $(S,q)$
 is finite; we shall use $\QD^1(S)$ to denote the space of \emph{unit-area} half-translation structures on $S$.

 Next, we consider geodesic representatives of a curve $\alpha$ on a half-translation surface $(S,q)$.
 In order to deal with punctured surfaces, we allow $\alpha$ to be homotoped so that it
 passes through punctures at the final moment of the homotopy, but not at any prior time.
 There are two possibilities: either $\alpha$ has a unique geodesic representative and is formed by concatenating
 a sequence of saddle connections (straight-line segments connecting singularities with no interior singularities),
 or there is a unique maximal (open) Euclidean cylinder foliated by the geodesic representatives of $\alpha$.
 In the latter case, we refer to $\alpha$ as a \emph{cylinder curve} on $(S,q)$.
 The \emph{Euclidean} (or \emph{flat}) length $l_q(\alpha)$ is the length of
 any geodesic representative of $\alpha$ on $(S,q)$ with respect to the Euclidean metric.
 The \emph{horizontal} and \emph{vertical} lengths of $\alpha$ on $(S,q)$,
 denoted by $l_q^H(\alpha)$ and $l_q^V(\alpha)$, are obtained by respectively
 integrating $|\Im(\sqrt{q})|$ and $|\Re(\sqrt{q})|$
 along a geodesic representative of $\alpha$.

 There is a natural $\SL$--action on $\QD^1(S)$ defined by $\R$--linear transformations of the natural co-ordinates. Two natural restrictions of this action yield the \emph{Teichm\"uller geodesic flow} and \emph{unipotent flow}, respectively defined by taking orbits under
 \[t \mapsto g_t := \begin{pmatrix}
  e^\frac{t}{2} & 0 \\
  0 & e^{-\frac{t}{2}}
  \end{pmatrix}
  \qquad \textrm{and} \qquad
  s \mapsto u_s := \begin{pmatrix}
  1 & s \\
  0 & 1
  \end{pmatrix};
 \]
 in particular, orbits under $g_t$ descend to (unit-speed) geodesics in $\T(S)$ called \emph{Teichm\"uller geodesics}.
 By Teichm\"uller's Theorem, there exists a unique Teichm\"uller geodesic connecting
 any given pair of distinct points in $\T(S)$.

 The $\SL$--orbit of quadratic differential $q\in\QD^1(S)$ descends to an isometrically embedded
 copy of the hyperbolic plane (of curvature $-4$) in $\T(S)$, called a \emph{Teichm\"uller disc},
 which we shall denote by $\Disc(q)$. The setwise stabiliser of $\Disc(q) \subset \T(S)$
 under the action of $\mcg(S)$ is naturally identified with the \emph{affine diffeomorphism group}
 $\Aff(q)$; an affine (self-)diffeomorphism of $(S,q)$ acts bijectively on the set of singularities,
 and restricts to a diffeomorphism with constant differential away from the singular set.
 The differential homomorphism $D \from \Aff(q) \to \PSL(2,\R)$ determines a short exact sequence 
 \[1 \to \Aut(q) \to \Aff(q) \to \PSL(q) \to 1, \]
 where the image $\PSL(q)$ is a Fuchsian group called the \emph{Veech group} \cite{Veech-dich}.
 The kernel is the the pointwise stabiliser of $\Disc(q)$, and coincides with the finite group
 of slope-preserving isometries of $(S,q)$.
 In the case where $\PSL(q)$ is a lattice then $(S,q)$ is called a \emph{Veech surface}. 
 Affine diffeomorphisms are classified as follows.

 \begin{proposition}[\cite{Thurston-diffeo, Veech-dich}]
  Let $\phi \in \Aff(q)$. Then
  \begin{enumerate}
   \item $D\phi$ is hyperbolic if and only if $\phi$ is pseudo-Anosov,
   \item $D\phi$ is parabolic if and only if some power of $\phi$ is a Dehn multitwist, and
   \item $D\phi$ is elliptic if and only if $\phi$ is periodic (finite order). \halmos
  \end{enumerate}
 \end{proposition}

 In the case where $D\phi$ is parabolic, there is an associated \emph{cylinder decomposition}
 of $(S,q)$. The foliation on $(S,q)$ by geodesics parallel to the unique eigenslope of $D\phi$
 is invariant under $\phi$; indeed, all separatrices with this slope are saddle connections.
 Cutting $(S,q)$ along these saddle connections decomposes it into a finite union of (maximal) Euclidean cylinders.
 The diffeomorphism $\phi$ acts by (possibly) permuting the cylinders, and by some power of a Dehn twist
 about the core curve of each cylinder.

 \subsection{Systoles and short markings}\label{sec:systole}

 Let us now review some natural maps from Teichm\"uller space to the curve complex and
 the marking complex.

 The \emph{extremal length} of $\alpha$ on $x \in \T(S)$ is
 \[\Ext_x(\alpha) := \sup_{\rho} \frac{l_\rho(\alpha)^2}{\area(\rho)},\]
 where $\rho$ runs over all metrics in the conformal class $x$.
 Given $\epsilon > 0$, we call $x \in \T(S)$ \emph{$\epsilon$--thick} if
 $\Ext_x(\alpha) \geq \epsilon$ for all curves $\alpha\in\C(S)$, and write
 $\T_{\geq \epsilon}(S)$ for the $\epsilon$--thick part of Teichm\"uller space.
 When $\epsilon$ is sufficiently small, $\T_{\geq \epsilon}(S)$ is connected.
 Moreover, by Mumford's Compactness Criterion \cite{Mumford-compactness, Bers-mumford}, $\mcg(S)$ acts geometrically on $\T_{\geq \epsilon}(S)$.

 The \emph{systole map} $\sys \from \T(S) \to \C(S)$ is defined by assigning $x\in\T(S)$
 its set of minimal length curves $\sys(x)\subset\C(S)$.
 The systole sets for $x\in\T(S)$ defined using extremal length, hyperbolic length,
 or Euclidean length for any $q \in \QD(x)$ agree in $\C(S)$ up to uniformly bounded Hausdorff distance (see \cite[Lemma 4.7]{TW-polygon}),
 and so, for our purposes, we may use whichever definition is most convenient.

 The \emph{short marking map} $\smark \from \T(S) \to \Mark(S)$ is defined as follows.
 The base curves of $\smark(x)$ are chosen to minimise hyperbolic length on $x\in\T(S)$
 according to the greedy algorithm. A transversal for each base curve $\beta$ is chosen
 to minimise $\dist_\beta(\cdot, a)$, where $a$ is a geodesic arc perpendicular to $\beta$
 with respect to the hyperbolic metric.
 The choice of short marking may not be unique, however, all the possible choices for
 $\smark(x)$ form a uniformly bounded diameter set in $\Mark(S)$.
 We shall also write $\mu_x$ to stand for $\smark(x)$.

 \begin{proposition}[\cite{MM1, MM2}]\label{prop:lipschitz}
  The maps $\sys \from \T(S) \to \C(S)$ and $\smark \from \T(S) \to \Mark(S)$
  are both coarsely Lipschitz. \halmos
 \end{proposition}

 Since $\smark$ is (coarsely) $\mcg(S)$--equivariant, it follows that
 the restriction $\smark \from \T_{\geq \epsilon}(S) \to \Mark(S)$ is a quasi-isometry.

 Subsurface projections can be defined on Teichm\"uller space by 
 composing with the short marking map.
 Note that $\pi_Y \circ \smark \from \T(S) \to \AC(Y)$ is also coarsely Lipschitz.
 As shorthand notation, write 
 \[\pi_Y(q) = \pi_Y(x) := \pi_Y(\mu_x)\]
 and
 \[\dist_Y(x,y) := \diam_{\AC(Y)}(\pi_Y(x) \cup \pi_Y(y))\]
 for $x, y \in \T(S)$ and $q\in\QD(x)$.

 \begin{theorem}[\cite{Rafi-teichdist}]\label{thm:rafi_dist}
  Fix a sufficiently small $\epsilon > 0$. Then there exists a constant $c_2 = c_2(S, \epsilon)$ such that
  for any $c \geq c_2$, there exists $\sA_2 > 0$ such that
  \[\dist_{\T(S)}(x,y) \asymp_{\sA_2} \sum_{Y \not\cong \Ann} \left[\dist_Y(x,y)\right]_c +
  \sum_{\alpha} \log \left[\dist_\alpha(x,y)\right]_c\]
  whenever $x,y \in \T(S)$ are $\epsilon$--thick.
  Here, the sums are respectively taken over all essential non-annular subsurfaces $Y \subseteq S$
  and all essential simple closed curves $\alpha$. \halmos
 \end{theorem}

 Let us now discuss the behaviour of the systole map along a Teichm\"uller geodesic.
 As shorthand notation, we shall write $t$ as subscript when referring to lengths on $q_t := g_t \cdot q$ (or the underlying conformal structure $x_t$).

 \begin{theorem}[\cite{Rafi-hyperbolicity}]\label{thm:shadow}
  Let $\G \from \R \to \T(S)$ be a Teichm\"uller geodesic and $Y\subseteq S$ be
  a subsurface. Then $\pi_Y \circ \G \from \R \to \AC(Y)$ is a uniform reparameterised
  quasigeodesic. Furthermore, if $Y \neq S$ then there exists a
  (possibly empty) interval $I_Y \subset \R$ such that
  \begin{itemize}
   \item each component of $~\R \setminus I_Y$ has uniformly bounded image under $\pi_Y \circ \G \from \R \to \AC(Y)$,
   and
   \item if $\alpha \subseteq \partial Y$ is a boundary component then $\Ext_t(\alpha) \leq \epsilon_0$ for all $t \in I_Y$. \halmos
  \end{itemize} 
 \end{theorem}

 Let $|I|$ denote the length of an interval $I \subseteq \R$.

 \begin{corollary}\label{cor:active}
 For every $L > 0$ there exists a constant $D = D(S,L) > 0$ such that the following holds.
 Suppose $Y$ is a proper subsurface that has a boundary component $\alpha$ satisfying
 $l_t(\alpha) \geq L$ for all $t \in \R$. Then
 $\dist_Y(\nu^+, \nu^-) \prec |I_Y| \prec D$.
 \end{corollary}

 \proof
 The condition $l_t(\alpha) \geq L$ for all $t \in \R$ implies that $\alpha$ cannot
 be completely horizontal nor completely vertical.
 The flat length of $\alpha$ satisfies
 \[l_t(\alpha) \asymp e^\frac{t}{2} l_q^H(\alpha) + e^{-\frac{t}{2}} l_q^V(\alpha) \asymp L_0\cosh(t - t_0)
 \geq L\cosh(t - t_0), \]
 where $t = t_0$ is a time at which $l_t(\alpha)$ attains a global minimum $L_0$.
 Using the definition of extremal length, we have
 \[ \Ext_t(\alpha) \succ L^2 \cosh^2(t - t_0).\]
 By the previous theorem, $\Ext_t(\alpha) \leq \epsilon_0$ for all $t \in I_Y$, and so
 $|I_Y|$ is bounded from above by some function of $\frac{\epsilon_0}{L^2}$.
 Combining this with the above theorem and the fact that $\pi_Y \from \T(S) \to \AC(Y)$
 is uniformly coarsely Lipschitz, we may bound $\dist_Y(\nu^+, \nu^-) \prec |I_Y|$ from above
 by some function of $\frac{\epsilon_0}{L^2}$.
 \endproof

 The \emph{no-backtracking property} follows from hyperbolicity of the curve complex and the fact that Teichm\"uller geodesics descend to (reparameterised) quasigeodesics in $\C(S)$.

 \begin{lemma}[\cite{MM1}]\label{lem:no_backtrack}
  There exists a constant $C = C(S)$ such that the following holds.
  Let $\G \from \R \to \T(S)$ be a Teichm\"uller geodesic.
  Then 
  \[\dist_S(\G(s), \G(u)) \geq \dist_S(\G(s), \G(t)) + \dist_S(\G(t), \G(u)) - C \]
  for all $s \leq t \leq u$. \halmos
 \end{lemma}

 In the case of thick Teichm\"uller geodesic segments, we have stronger control over their images in $\C(S)$ using Theorem \ref{thm:rafi_dist} and Rafi's characterisation of short curves
 along Teichm\"uller geodesics.

 \begin{proposition}[\cite{Rafi-short, Rafi-teichdist}]\label{prop:thick_geod}
  Let $\epsilon > 0$ be sufficiently small, and suppose $\G \from I \to \T(S)$ is an $\epsilon$--thick
  Teichm\"uller geodesic segment.
  Then ${\sys \circ \G \from I \to \C(S)}$ is a $\sK$--quasigeodesic, where $\sK = \sK(S,\epsilon)$.
  \halmos
 \end{proposition}

 We finish off the section by briefly recalling the notion of the \emph{geodesic representative}
 $\sY_q$ of a subsurface $Y \subseteq S$ on a half-translation structure $(S,q)$ due to Rafi \cite{Rafi-short} (see also \cite{MT-fibered}).
 First suppose $Y$ is an annulus. If the core curve of $Y$ is a cylinder curve, then let $\sY_q$ be
 the maximal Euclidean cylinder foliated by geodesic representatives of the core curve;
 otherwise, set $\sY_q$ to be the unique geodesic representative of its core curve.
 Now suppose $Y$ is non-annular. The idea is to ``pull tight'' each boundary component of $\partial Y$ to
 their geodesic representatives to obtain $\sY_q$.
 If a component of $\partial Y$ is a cylinder curve, then we require $\sY_q$
 to be disjoint from the interior of the associated maximal cylinder.
 If some component $\alpha$ of $\partial Y$ is homotopic to a puncture of $S$, then the geodesic representative
 of $\alpha$ to degenerates to the completion point on $(S,q)$ associated to the puncture.
 Following Minsky--Taylor, we say that $Y$ is \emph{$q$--compatible} if there is a homotopy between $Y$ and $\sY_q$ that restricts to
 an isotopy between their interiors. If this holds, then cutting $(S,q)$ along all saddle connections appearing on
 $\partial \sY_q$ yields $\sY_q$ as a complementary component.

 \begin{proposition}[\cite{MT-fibered}]\label{prop:q-compatible}
  Let $\nu^+, \nu^-$ be the horizontal and vertical foliations associated to $q \in \QD^1(S)$. If a subsurface $Y$ is not $q$--compatible then $\dist_Y(\nu^+, \nu^-) \leq 3$ (and the subsurface projections $\pi_Y(\nu^\pm)$ are non-empty).\halmos
 \end{proposition}

 \section{Bounded geometry and Nielsen cores}

 For the rest of this paper, we shall fix a finitely generated subgroup $\hat\Gamma \leq \mcg(S)$
 that stabilises a Teichm\"uller disc $\Disc(\Gamma) \subset \T(S)$.
 Furthermore, we assume that $\hat\Gamma$ is not virtually cyclic, as such groups are known to
 be undistorted in $\mcg(S)$ \cite{FLM-rank1}.
 This also guarantees that the Teichm\"uller disc is unique; see the remark below.

 The group $\hat\Gamma$ can also be viewed as a subgroup of $\Aff(q)$ for any $q\in\QD(S)$
 generating $\Disc(\Gamma)$. 
 Restricting the differential homomorphism to $\hat\Gamma$ gives rise to a short exact equence
 \[1 \to \hat\Gamma \cap \Aut(q) \to \hat\Gamma \to \Gamma \to 1,\]
 where the image $\Gamma \leq \PSL(q)$ is a finitely generated non-elementary Fuchsian group.
 Since the kernel is finite, the quotient map $\hat \Gamma \to \Gamma$ is a quasi-isometry
 with respect to their word metrics. Therefore, for our purposes, we may equally work with
 the actions of $\hat \Gamma$ or $\Gamma$ on the associated Teichm\"uller disc $\Disc(\Gamma)$,
 regarded respectively as either a subset of $\T(S)$ or as a copy of $\Hy^2$.

 All quadratic differentials $q$ shall henceforth be assumed to belong to the
 $\SL$--orbit in $\QD^1(S)$ descending to $\Disc(\Gamma)$
 unless otherwise specified.

 \begin{remark}
 Any non-elementary Fuchsian group contains a hyperbolic element, and so $\hat\Gamma$
 contains a pseudo-Anosov element.
 Any pseudo-Anosov element stabilises a unique Teichm\"uller geodesic, and thus determines a
 unique Teichm\"uller disc. Therefore, the Teichm\"uller disc $\Disc(\Gamma)$ is uniquely
 determined by the subgroup $\Gamma$.
 \end{remark}

 \begin{definition}
 The \emph{Nielsen core} $\Core(\Gamma) \subseteq \Disc(\Gamma)$ of $\Gamma$ is the convex hull of
 the limit set ${\Lambda(\Gamma) \subseteq \partial \Disc(\Gamma)}$.
 \end{definition}

 The inclusion $\Core(\Gamma) \hookrightarrow \Disc(\Gamma)$ is an isometric embedding.
 Since the action of $\mcg(S)$ on $\T(S)$
 is properly discontinuous, the same holds for the action of $\Gamma$ on $\Core(\Gamma)$.
 The quotient $\Core(\Gamma)/\Gamma$ is a finite area hyperbolic orbifold, possibly with
 geodesic boundary.
 Since $\Gamma$ is finitely generated, $\Core(\Gamma)/\Gamma$ has empty boundary precisely when 
 $\Lambda(\Gamma) = \partial\Disc(\Gamma)$, in which case $\Gamma$ is a lattice.

 The goal of this section is to establish positive constants  which control the geometry of quadratic differentials appearing over $\Core(\Gamma)$.

 \subsection{Parabolic subgroups}

 The (finite) set of cusps of $\Core(\Gamma)/\Gamma$ is in one-to-one correspondence with the
 conjugacy classes of maximal parabolic subgroups of $\Gamma$.
 Let $\Par(\Gamma)$ be the set of all maximal parabolic subgroups of
 $\Gamma$. For each $H \in \Par(\Gamma)$, let $\ParC(H)$ be the set of core curves of the
 cylinders associated to $H$. 
 Define $\ParC(\Gamma) := \cup_{H \in \Par(\Gamma)} \ParC(H)$ to be the set of \emph{parabolic cylinder
 curves} associated to $\Gamma$. Furthermore, if a saddle connection is parallel to some parabolic
 slope then we shall call it a \emph{parabolic saddle connection}.

 \begin{remark}
  If $\Gamma$ has at least one (maximal) parabolic subgroup $H$
  then it must have infinitely many; these can be obtained, for example, by conjugating $H$ by powers of a hyperbolic
  element of $\Gamma$. It may be the case that $\Par(\Gamma)$ is empty, occuring precisely when $\Core(\Gamma)/\Gamma$ has no cusps;
  in this situation $\Core(\Gamma)/\Gamma$ has at least one geodesic boundary component as $\Gamma$ cannot act cocompactly on $\Disc(\Gamma)$.
 \end{remark}

 \subsection{Cylinder widths}

 Given a parabolic subgroup $H \in \Par(\Gamma)$,
 let $W_H(q) > 0$ be the minimum width of all (maximal) cylinders on $(S,q)$ whose core curve belongs to $\ParC(H)$.
 Since $(S,q)$ has unit area, we deduce that each curve in $\ParC(H)$ has flat length at most $\frac{1}{W_H(q)}$ on $(S,q)$.  Note that $W_H(q) \rightarrow \infty$ as $q$ tends towards the fixed point of $H$ on $\partial \Disc(\Gamma)$.
 It follows that the function 
 \[q \mapsto \sup_{H \in \Par(\Gamma)} W_H(q)\]
 defined on $\Core(\Gamma)$
 descends to a continuous proper function on $\Core(\Gamma)/\Gamma$,
 and thus attains a positive minimum value $W_\Gamma > 0$ at some $q_0\in \Core(\Gamma)$.
 Now, there are only finitely many curves on $(S,q_0)$ whose flat length is bounded above
 by any given positive constant. Therefore, there are finitely many $H \in \Par(\Gamma)$
 for which $W_H(q_0)$ is bounded from below by any given positive constant,
 and so the supremum is attained by some parabolic subgroup.

 \begin{lemma}\label{lem:non_par_slope}
 Every non-parabolic saddle connection on $q \in \Core(\Gamma)$ has length at least $W_\Gamma$.
 In particular, if $\alpha\in\C(S)$ is not parallel to a parabolic slope, then $l_q(\alpha) \geq W_\Gamma$ throughout $\Core(\Gamma)$.
 \end{lemma}

 \proof
 For each $q \in \Core(\Gamma)$, there exists a parabolic subgroup $H \in \Par(\Gamma)$ for which
 every cylinder in $\ParC(H)$ has width at least $W_\Gamma$ on $(S,q)$.
 Then any non-parabolic saddle connection on $(S,q)$ must intersect some cylinder in $\ParC(H)$ transversely,
 and thus has length at least $W_\Gamma$.
 If $\alpha$ is not parallel to a parabolic slope, then its geodesic representative on $(S,q)$ must
 use at least one saddle connection not parallel to the slope of $H$. The result follows.
 \endproof

 \begin{corollary}\label{cor:ext_width}
  If $\alpha \in \C(S)$ is not parallel to a parabolic slope then $\Ext_q(\alpha) \geq W_\Gamma^2$
  throughout $\Core(\Gamma)$. \halmos
 \end{corollary}

 \subsection{Expanding annuli}

 Let us turn our attention to curves that are parallel to some parabolic slope for $\Gamma$.
 These curves can have arbitrarily short flat length on $\Core(\Gamma)$.
 Our goal is to show that when such a curve is not itself a parabolic cylinder curve
 then its extremal length is bounded below by some constant depending only on $\Gamma$.

 We shall briefly recall the notions of flat and expanding annuli from Minsky \cite{Minsky-harmonic}.
 Given a curve $\alpha\in\C(S)$, define its (possibly degenerate) \emph{flat annulus} $F_q(\alpha)$ as follows:
 If $\alpha$ is a cylinder curve on $(S,q)$, then let $F_q(\alpha)$ be the associated maximal
 flat cylinder; otherwise let $F_q(\alpha)$ be the (unique) geodesic representative of $\alpha$
 on $(S,q)$. Note that $F_q(\alpha)$ contains all geodesic representatives of $\alpha$.

 Next, equip the annular cover $S^\alpha$ with the pullback metric from $(S,q)$.
 The flat annulus $F_q(\alpha)$ lifts to a unique flat annulus $\tilde F_q(\alpha)$
 on $S^\alpha$. Cutting $S^\alpha$ along the two (possibly coincident) boundary curves of
 $\tilde F_q(\alpha)$ yields exactly two components $S^\alpha_+$ and $S^\alpha_-$ that are not flat cylinders.
 Given $r \geq 0$, let $\tilde E^\pm_q(r) \subset S^\alpha_\pm$ be the intersection of $S^\alpha_\pm$
 with the closed {$r$--neighbourhood} of $\tilde F_q(\alpha)$ in $S^\alpha$.
 For $r > 0$, $\tilde E^\pm_q(r)$ is topologically a closed annulus.

 Now, consider the projection of $\tilde E^\pm_q(r)$ to $(S,q)$.
 If the projection map is injective on the interior of $\tilde E^\pm_q(r)$, then
 we call the image $E^\pm_q(r)$ a \emph{(regular) expanding annulus} of $\alpha$ on $(S,q)$ of \emph{radius} $r$.
 Furthermore, if $E^\pm_q(r)$ has no singularities in its interior then it is called \emph{primitive};
 let $r_q^\pm(\alpha)\geq 0$ be the largest value of $r$ for which this holds. 
 Call the boundary curve of $E^\pm_q(r)$ that coincides with a boundary curve of $F_q(\alpha)$ the \emph{inner boundary},
 and the other boundary curve the \emph{outer boundary} (these will coincide when $r=0$).

 The following theorem gives an estimate for the extremal length of extremely short curves in terms of the geometry
 of their maximal flat and primitive expanding annuli.
 Recall that the \emph{modulus} $Mod(C)$ of a flat cylinder $C$ is its width divided by the length of its core curve.

 \begin{theorem}[\cite{Minsky-harmonic, CRS-minima}]\label{thm:CRS}
 If $\alpha$ is extremely short on $(S,q)$ then   \[\frac{1}{\Ext_q(\alpha)} \asymp \max\left\{Mod(F_q(\alpha)), ~\log \left(\frac{r^\pm_q(\alpha)}{l_q(\alpha)}\right) \right\}. \tag*{\qed}\]
 \end{theorem}

 Let us now focus our attention on primitive expanding annuli associated to curves parallel to some parabolic slope for $\Gamma$.
 Given a parabolic subgroup $H \in \Par(\Gamma)$,
 let $\gamma_H, \gamma_H'$ respectively be a shortest and longest saddle connections on $(S,q)$ parallel
 to the slope corresponding to $H$, for some $q \in \Disc(\Gamma)$.
 (These saddle connections will respectively remain shortest and longest throughout $\Disc(\Gamma)$.)
 Define 
 \[\rho_H := \frac{l_q(\gamma_H')}{l_q(\gamma_H)} > 0,\]
 for any $q \in \Disc(\Gamma)$.
 Note that this ratio is constant under $\SL$--deformations.
 Since there are finitely many conjugacy classes of maximal parabolic subgroups,
 it follows that the supremum
 \[\rho_\Gamma := \sup_{H \in \Par(\Gamma)} \rho_H < \infty.\]
 is finite and attained.

 \begin{lemma}
  Let $\alpha$ be a curve that is parallel to some parabolic slope for $\Gamma$.
  Then 
  \[\frac{r^\pm_q(\alpha)}{l_q(\alpha)} \leq \rho_\Gamma\]
  for all $q \in \Core(\Gamma)$.
 \end{lemma}

 \proof
 By applying a rotation, we may assume that $\alpha$ is horizontal on $(S,q)$.
 Since the horizontal slope is parabolic, $(S,q)$ admits a horizontal cylinder decomposition.
 In particular, every horizontal separatrix is a saddle connection.

 Let $E$ be a maximal primitive expanding annulus for $\alpha$ on $(S,q)$.
 If $r^\pm_q(\alpha) = 0$ then we are done, so we may assume that $E$ has non-empty interior.
 The inner boundary of $E$ is a geodesic representative of $\alpha$ passing through at least one singularity.
 Since $\alpha$ is horizontal, the interior angle (inside $E$) at each singularity on the inner boundary
 is an integer multiple of $\pi$.
 We claim that at least one such singularity has interior angle at least $2\pi$.
 If not, then, for sufficiently small $0 < r < r^\pm_q(\alpha)$, the annulus $E^\pm_q(r) \subset E$ is isometric
 to a Euclidean cylinder; this contradicts the construction of flat and expanding annuli.
 Therefore, there exists a horizontal saddle connection $\beta$ starting on the inner boundary of $E$
 with an initial segment lying in the interior of $E$. 
 Since $E$ is primitive, we deduce that $r_q^\pm(\alpha) \leq l_q(\beta)$, for otherwise $E$
 will contain an interior singularity.
 On the other hand, $l_q(\alpha)$ is at least the length of the shortest horizontal saddle connection.
 The desired result follows using the definition of $\rho_\Gamma$.
 \endproof

 Any curve parallel to a parabolic slope that is not itself a parabolic cylinder curve has a degenerate flat annulus.
 Thus, combining the preceding results with Corollary \ref{cor:ext_width} yields the following.

 \begin{proposition}\label{prop:short}
  There exists a constant $0 < \epsilon_\Gamma \leq \epsilon_0$ such that for every curve $\alpha \in \C(S)$,
  either 
  \begin{itemize}
  \item $\alpha$ is a parabolic cylinder curve for $\Gamma$ and $\inf_{q\in\Core(\Gamma)} \Ext_q(\alpha) = 0$, or
  \item $\alpha$ is not a parabolic cylinder curve for $\Gamma$ and $\Ext_q(\alpha) \geq \epsilon_\Gamma$ for all $q \in \Core(\Gamma)$. \halmos
  \end{itemize}
 \end{proposition}

 In the following subsections, we may need to take $\epsilon_\Gamma$ sufficiently small
 to ensure that desired properties hold. 

 \subsection{Horodiscs}

 We now choose a preferred family of horodiscs in $\Disc(\Gamma)$ associated to the family
 of parabolic subgroups.
 For each $H \in \Par(\Gamma)$, 
 let $\gamma_H$ be a shortest saddle connection with slope corresponding to $H$,
 then define a pair of nested horodiscs
 \[U(H) = \set{q \in \Disc(\Gamma) \st l_q(\gamma_H) \leq \sqrt{\epsilon_\Gamma}} \quad \textrm{and} \quad
 U'(H) = \set{q \in \Disc(\Gamma) \st l_q(\gamma_H) \leq \sqrt{\epsilon_0}}.\] 
 Each $U(H)$ descends to a neighbourhood of a cusp on $\Core(\Gamma)/\Gamma$.
 Since there are finitely many cusps, we may choose $\epsilon_\Gamma$ sufficiently
 small to ensure that 
 \[\dist_{\Disc(\Gamma)}(U(H), U(K)) \geq 1\]
 for all distinct $H,K \in \Par(\Gamma)$.
 This also ensures that the cusp neighbourhood arising from each $U(H)$ is topologically an annulus.
 Next, define the \emph{truncated Nielsen core} of $\Gamma$ to be
 \[\trCore(\Gamma) := \Core(\Gamma) \setminus \bigcup_{H\in\Par(\Gamma)} int(U(H)).\]
 By Proposition \ref{prop:short}, the $\epsilon_\Gamma$--thin part of $\Core(\Gamma)$
 is contained in $\bigcup_{H\in\Par(\Gamma)} U(H)$,
 and so $\trCore(\Gamma)$ is $\epsilon_\Gamma$--thick.

 \begin{lemma}\label{lem:horo_overlap}
  There exists a constant $R_\Gamma > 0$ such that for all distinct $H,K \in \Par(\Gamma)$,
  we have $\diam(U'(H) \cap U'(K)) \leq R_\Gamma$.
 \end{lemma}

 \proof
 Let $t = \log(\epsilon_0/\epsilon_\Gamma) \geq 0$.
 Suppose $q \in U'(H)$. By applying a rotation, we may assume
 $\gamma_H$ is vertical on $(S,q)$.
 Then 
 \[l_{g_t\cdot q}(\gamma_H) = e^{-t/2}l_q(\gamma_H) \leq \sqrt{\epsilon_\Gamma/\epsilon_0}\sqrt{\epsilon_0} = \sqrt{\epsilon_\Gamma}.\]
 Therefore, $U'(H)$ is contained in the $t$--neighbourhood of $U(H)$
 for all $H \in \Par(\Gamma)$.
 The diameter bound then follows using elementary hyperbolic geometry and the fact that
 the distance between distinct horodiscs $U(H)$ and $U(K)$ is at least 1.
 \endproof

 \subsection{Virtual triangle areas}

 Smillie and Weiss characterise Veech surfaces as the half-translation surfaces 
 whose \emph{virtual triangle spectrum} is discrete \cite{SW-Veech}.
 Motivated by their work, we consider the \emph{parabolic virtual triangle spectrum}
 of $\Gamma$ obtained by restricting to the parabolic saddle connections, and prove that it is always discrete.
 This is not necessary for our main theorem, however, we include
 it as it may be of independent interest.

 Associated to any saddle connection on $(S,q)$ is a \emph{holonomy vector} in $\CC$
 that has the same slope and length; this is well-defined up to scaling by $\pm 1$.
 Let $\hol_\Gamma(q)$ be the set of holonomy vectors associated to parabolic saddle
 connections on $(S,q)$ (with respect to $\Gamma$).
 Define the \emph{parabolic virtual triangle spectrum} of~ $\Gamma$ to be
 \[PVT(\Gamma) := \set{|u \wedge v| \st u, v \in \hol_\Gamma(q) } \subset \R,\]
 for some (hence all) $q \in \Disc(\Gamma)$.
 Let $v_H(q)$ be the holonomy vector of a shortest saddle connection parallel to the
 slope corresponding to $H \in \Par(\Gamma)$.

 \begin{lemma}
  Let $H, K \in \Par(\Gamma)$. Then
  \[\dist_{\Disc(\Gamma)}(U(H), U(K)) = 2\log\left( \frac{|v_H \wedge v_K|}{\epsilon_\Gamma} \right). \]
 \end{lemma}

 \proof
 Let $\G(t) = q_t$ be the infinite Teichm\"uller geodesic whose horizontal and vertical
 slopes correspond to the slopes of $H$ and $K$ respectively, and where $q_0 \in \partial U(H)$.
 Note that the unique geodesic segment connecting $U(H)$ and $U(K)$ is a subinterval of $\G$.
 Using the definition of the horodiscs, we have $v_H(q_0) = \sqrt{\epsilon_\Gamma}$
 and $v_K(q_0) = e^{\frac{t}{2}} v_K(q_t) = e^{\frac{t}{2}} \sqrt{\epsilon_\Gamma}$
 when $t = \dist_{\Disc(\Gamma)}(U(H), U(K))$. The result follows.
 \endproof

 \begin{lemma}
  The set $PVT(\Gamma)$ is discrete in $\R$.
 \end{lemma}

 \proof
 It suffices to show that $PVT(\Gamma) \cap [0, b]$ contains finitely many values for all $b \geq 0$.
 By the above lemma, this is equivalent to proving that there are
 finitely many geodesic segments on $\trCore(\Gamma)/\Gamma$
 orthogonal to cusp boundaries of length less than any given positive constant.
 This follows from a standard argument; for example, by doubling $\trCore(\Gamma)/\Gamma$ along its boundary
 to obtain a closed compact surface, and using the fact that there are finitely many closed curves whose
 geodesic length is bounded above by any given constant.
 \endproof

 \subsection{Bounded projection image}\label{sec:bounded}

 The main technical result of this paper is  the following dichotomy for subsurface projections
 Call an annulus on $S$ \emph{parabolic} for $\Gamma$ if its core curve is a parabolic cylinder curve for $\Gamma$.

 \newpage
 \begin{proposition}\label{prop:diameter}
 There exists a constant $D_\Gamma > 0$ such that given any subsurface  $Y \subsetneq S$, the set ${\pi_Y( \smark(\Core(\Gamma)) )}$ has
 \begin{itemize}
  \item infinite diameter in $\AC(Y)$ if $Y$ is a parabolic annulus for $\Gamma$, and
  \item diameter at most $D_\Gamma$ in $\AC(Y)$ otherwise.
 \end{itemize}
 \end{proposition}

 For brevity, write $M(\Gamma) := \smark(\Core(\Gamma))$.

 First, we consider the case where $Y$ is a parabolic annulus.
 Let $H$ be the parabolic subgroup containing $Y$ in its associated cylinder decomposition.
 Let $\phi\in H$ be a non-trivial element that preserves $Y$.
 Note that $\phi|_Y$ is a power of a Dehn twist about the core curve of $Y$.
 Furthermore, $\pi_Y(\phi^n \cdot \mu) = (\phi|_Y)^n \cdot \pi_Y(\mu)$ for any marking $\mu \in M(\Gamma)$.
 Since orbits in $\AC(Y)$ under $\langle \phi|_Y \rangle$ have infinite diameter, it follows that
 $\pi_Y(M(\Gamma)) \supseteq \pi_Y(\langle \phi \rangle \cdot \mu)$ is unbounded.

 We may henceforth assume $Y \neq S$ is not a parabolic annulus.
 By Proposition \ref{prop:q-compatible}, we may also assume that $Y$ is
 $q$--compatible.

 Let $\bG(\Gamma)$ be the set of bi-infinite geodesics on $\Disc(\Gamma)$ with
 both endpoints in the limit set $\Lambda(\Gamma)$. 
 Note that $\Core(\Gamma)$
 is the convex hull of $\bigcup_{\G \in \bG(\Gamma)} \G$ in $\Disc(\Gamma)$.
 For concreteness, we take the images of $U(H)$ on $\Core(\Gamma)/\Gamma$
 to be a set of preferred cusp neighbourhoods, where $H$ runs over a set of representives
 for each conjugacy class in $\Par(\Gamma)$.

 \begin{lemma}
  There exists some $\sr_\Gamma > 0$ such that $\Core(\Gamma)$ is contained
  in the $\sr_\Gamma$--neighbourhood of ~$\bigcup_{\G \in \bG(\Gamma)} \G$ in $\Disc(\Gamma)$.
 \end{lemma}

 \proof
 It suffices to prove the result in the case where $\Disc(\Gamma)/\Gamma$ is a hyperbolic surface;
 the general case can be dealt with by taking a finite orbifold cover.
 In this situation, every geodesic $\G \in \bG(\Gamma)$ descends to a complete geodesic on $\Core(\Gamma)/\Gamma$;
 moreover, every complete geodesic on $\Core(\Gamma)/\Gamma$ arises this way.
 Note that a complete geodesic on $\Core(\Gamma)/\Gamma$ cannot be contained in any cusp neighbourhood,
 and so must have non-empty intersection with $\trCore(\Gamma)/\Gamma$.
 It follows that $\trCore(\Gamma)/\Gamma$ is contained in the $\sr$--neighbourhood of any complete geodesic
 on $\Core(\Gamma)/\Gamma$ for $\sr \geq \diam(\trCore(\Gamma)/\Gamma)$.

 We now deal with the cusp neighbourhoods of $\Core(\Gamma)/\Gamma$.
 Since $\Gamma$ is non-elementary, it either has no (maximal) parabolic subgroups, or infinitely many.
 Therefore, for every $H \in \Par(\Gamma)$ there exists some geodesic $\G \in \bG(\Gamma)$
 with one end contained in $U(H)$. The associated cusp neighbourhood in $\Core(\Gamma)/\Gamma$ is then contained in the
 $\sr$--neighbourhood of the image of $\G$, as long as $\sr > 0$ is greater than the length of the cusp boundary.

 The desired result holds for any value of $\sr_\Gamma > 0$ greater than both the diameter of $\trCore(\Gamma)/ \Gamma$,
 and its longest cusp boundary length.
 \endproof

 Using Proposition \ref{prop:lipschitz}, the following is immediate.

 \begin{corollary}
  There exists some $\sr'_\Gamma > 0$ such that for all $Y$, the set $\pi_Y (M(\Gamma))$ is contained in
  the $\sr'_\Gamma$--neighbourhood of ~$\bigcup_{\G \in \bG(\Gamma)} \pi_Y(\G)$ in $\AC(Y)$. \halmos
 \end{corollary}

 Our strategy is to now show that $\pi_Y(\Lambda(\Gamma))$ has diameter in $\AC(Y)$ bounded above
 by a constant independent of the choice of $Y$.
 Then for any $\G \in \bG(\Gamma)$ with endpoints $\nu^+$ and $\nu^-$, the uniform reparameterised quasigeodesic $\pi_Y \circ \G$ (coarsely) connects $\pi_Y(\nu^+)$ and $\pi_Y(\nu^-)$ in $\AC(Y)$ and hence has uniformly bounded diameter (see Theorem \ref{thm:shadow}).
 Observe that for any pair $\G, \G' \in \bG(\Gamma)$, there exists some $\G'' \in \bG(\Gamma)$
 sharing at least one endpoint with each of $\G$ and $\G'$; this implies that $\bigcup_{\G \in \bG(\Gamma)} \pi_Y(\G)$ has diameter bounded from above by a constant depending only on $\Gamma$.
 Appealing to the above corollary completes the proof of Proposition \ref{prop:diameter}.

 \begin{proposition}
  There exists some $D > 0$ such that for any subsurface $Y \neq S$, not a parabolic annulus,
  the set $\pi_Y(\Lambda(\Gamma))$ has diameter at most $D$ in $\AC(Y)$.
 \end{proposition}

 \proof
 Let $Y \neq S$ be a subsurface that is not a parabolic annulus.
 The proof proceeds in three cases depending on the geodesic representatives of $\partial Y$.

 \subsubsection*{Case 1: $\partial Y$ has a boundary component $\gamma$ that is not parallel to some parabolic
 slope.} 

 Let $\nu^+, \nu^- \in \Lambda(\Gamma)$ be a pair of distinct foliations,
 and $\G \in \bG(\Gamma)$ be the Teichm\"uller geodesic in $\Core(\Gamma)$ connecting them.
 By Lemma \ref{lem:non_par_slope}, we have $l_t(\gamma) \geq W_\Gamma$ for all $t$.
 Applying Corollary \ref{cor:active}, we deduce that
 $\dist_Y(\nu^+, \nu^-) \prec |I_Y| \prec D_1$
 for some $D_1 = D_1(S,W_\Gamma)$.

 \subsubsection*{Case 2: Each boundary curve of $Y$ has a parabolic slope, but they are not all parallel.}

 Let $\gamma_1, \gamma_2 \subset \partial Y$ be boundary curves parallel to distinct parabolic slopes,
 and let $H_1, H_2 \in \Par(\Gamma)$ respectively be the corresponding parabolic subgroups.
 By Theorem \ref{thm:shadow} and the definition of $U'(H_i)$,
 we have ${\G(I_Y) \subseteq U'(H_1) \cap U'(H_2)}$ for any $\G\in\bG(\Gamma)$.  Then $|I_Y| \leq R_\Gamma$, by Lemma \ref{lem:horo_overlap},
 and so the result follows using Corollary \ref{cor:active}.

 \subsubsection*{Case 3: All boundary curves of $Y$ are parallel and have parabolic slope.}

 For the remaining case, note that $Y$ cannot be an annulus.
 We shall prove a stronger statement in order to bound $\dist_Y(\nu^+, \nu^-)$.
 Given $q\in\QD(S)$, let $\PMF(q)$ be the set of projectivised measured foliations
 arising as the horizontal foliation of $e^{i\theta} q$ for some $\theta\in\RP$.
 Recall that $\sY_q$ is the geodesic representative of $Y$ on $(S,q)$.

 \begin{lemma}
 Let $q\in \QD(S)$ and suppose $Y \neq S$ is a non-annular $q$--compatible subsurface with horizontal boundary.
 Let $\gamma \in \C(Y)$ be a horizontal curve. Then
 \[\diam_Y \PMF(q) \prec \log\left(\frac{l_q(\gamma)}{l_q(\partial Y)}\right).\]
 \end{lemma}

 Note, this bound is silent if there exist no essential horizontal curves on $Y$.

 \proof 
 The strategy is to bound $\dist_Y(\gamma, \nu)$ from above for all $\nu \in \PMF(q)$.   If $\nu$ is horizontal then $\dist_Y(\gamma, \nu) \leq 1$, so we may assume otherwise.
 By applying an appropriate $\SL$--deformation, we may arrange so that
 $\nu$ is vertical while preserving the horizontal slope.
 Cut $\sY_q$ along all vertical separatrices that start either at a boundary singularity
 with internal angle at least $2\pi$ or an interior singularity, and end either at a singularity or on the boundary of $\sY_q$.
 This decomposes $\sY_q$ into a union of (at least one) Euclidean rectangles (with horizontal and vertical sides),
 and a (possibly empty) set of subsurfaces with vertical boundary.
 The number of separatrices that were cut along is bounded above in terms of $|\chi(Y)| \leq |\chi(S)|$,
 and so the number of rectangles is also bounded above in terms of $|\chi(Y)|$.
 Note that $\partial \sY_q$ is the union of the horizontal sides of these rectangles, and so
 the widths of these rectangles sum to $2l_q(\partial Y)$.
 Therefore, there exists a rectangle $R$ of $\width(R) \succ l_q(\partial Y)$. Let $\eta$ be the vertical arc in $R$ that connects the midpoints of
 its two horizontal sides. 
 Then 
 \[i(\eta, \gamma) \leq \frac{l_q(\gamma)}{\width(R)} \prec \frac{l_q(\gamma)}{l_q(\partial Y)},\]
 and so by Hempel's bound (\ref{eqn:hempel}) we have
 \[\dist_Y(\eta, \gamma) \prec \log \left(\frac{l_q(\gamma)}{l_q(\partial Y)}\right).\]
 Since $\nu$ is vertical, it has no transverse intersection with $\eta$, and so
 $\dist_Y(\eta,\nu) \leq 1$.
 \endproof

 It remains to show that $Y$ has an essential curve that is not too long
 compared to the length of $\partial Y$.
 By assumption, $Y$ is $q$--compatible and thus has embedded interior.
 Since the horizontal slope is parabolic, cutting $(S,q)$ along all horizontal
 saddle connections decomposes it into a union of horizontal cylinders.
 In particular, $\sY_q$ can be formed by taking a non-empty subset of these cylinders,
 then gluing them along some horizontal saddle connections; let $\gamma$ be
 a shortest core curve of a horizontal cylinder contained in $\sY_q$.
 Note that $\gamma$ cannot be peripheral on $Y$, for if any boundary component of $Y$
 is a cylinder curve, then the interior of the associated maximal cylinder must be disjoint
 from $\sY_q$.

 We wish to bound 
 $\frac{l_q(\gamma)}{l_q(\partial Y)}$ from above.
 Let $l_0$ and $l_1$ respectively be the lengths of the shortest and the longest horizontal
 saddle connection on $(S,q)$.
 Note that $l_q(\partial Y) \geq l_0$. Consider the maximal cylinder $C$ with core curve $\gamma$.
 Each boundary component of $C$ runs over any given saddle connection on $(S,q)$ at most twice.
 Since the number of horizontal saddle connections on $(S,q)$ is bounded above in terms of $|\chi(S)|$,
 we deduce that $l_q(\gamma) \prec l_1$.
 Therefore 
 \[ \frac{l_q(\gamma)}{l_q(\partial Y)} \prec \frac{l_1}{l_0} \leq\rho_\Gamma.\]
 Applying the above lemma completes the proof of Proposition \ref{prop:diameter}. \halmos

 \subsection{Cusp winding}

 In this section, we estimate the annular projection distance
 $\dist_\alpha(q,q')$ for points $q,q' \in \trCore(\Gamma)$ and $\alpha\in\ParC(\Gamma)$
 in terms of the amount of winding about the associated cusp.

 Let us recall Rafi's estimate for annular projection distance in terms of the \emph{relative twisting}
 of a pair of quadratic differentials about a curve $\alpha \in \C(S)$.
 Given $q \in \QD(S)$, let $\eta_\alpha(q)$ be a complete Euclidean geodesic on $(S,q)$ orthogonal
 to the geodesic representative of $\alpha$.
 (We may also assume that $\eta_\alpha(q)$ does not hit any singularities.)
 Let $\tilde\eta_\alpha(q)$ be a lift of $\eta_\alpha(q)$ on the annular
 cover $S^\alpha$ that intersects the unique closed lift of $\alpha$ essentially.
 Then the relative twisting $\tw_\alpha(q, q')$ is defined to be the geometric intersection number between
 $\tilde\eta_\alpha(q)$ and $\tilde\eta_\alpha(q')$; this is well-defined up to a uniform additive error.

 \begin{proposition}[\cite{Rafi-hyperbolicity}]\label{prop:rafi_twist}
  For all $q,q' \in \QD(S)$ and $\alpha \in \C(S)$ we have $\tw_\alpha(q, q') \asymp \dist_\alpha(q, q')$. \halmos
 \end{proposition}

 We now focus on the case where $\alpha$ is a parabolic cylinder curve for $\Gamma$,
 and where the quadratic differentials are restricted to $\Disc(\Gamma)$.
 Choose $H \in \Par(\Gamma)$ so that ${\alpha \in \ParC(H)}$,
 and suppose $\G$ is a Teichm\"uller geodesic on $\Disc(\Gamma)$
 orthogonal to $\partial U(H)$. By applying a suitable rotation, we may assume that
 all cylinders in $\ParC(H)$ are horizontal on $q_t \in \G$ for all $t \in \R$.
 Therefore, any geodesic on $(S,q_t)$ orthogonal to the geodesic representative of $\alpha$
 is vertical, and hence is a leaf of the vertical foliation on $(S,q_t)$.
 In particular, $\eta_\alpha(q_t)$ can be chosen to be the same topological leaf for all $t\in\R$. 
 It follows that $\tw_\alpha(q, q') = 0$ for all $q,q'\in\G$.

 Next, we define the \emph{cusp winding} with respect to a parabolic subgroup $H \in \ParC(\Gamma)$ as follows:
 Given $q, q' \in \Disc(\Gamma)$, let $\dist_H(q, q')$ be the distance between 
 their respective nearest point projections to the horocycle $\partial U(H)$, measured along $\partial U(H)$. 
 Note that if $\G$ is a geodesic in $\Disc(\Gamma)$ orthogonal to $\partial U(H)$
 then all points along $\G$ project to a common point on $\partial U(H)$;
 thus $\dist_H$ gives a notion of distance between pairs of such geodesics.

 We shall show that the relative twisting and cusp winding agree up to uniform additive and multiplicative
 error depending only on $\Gamma$. Let us introduce some more constants.
 Given a parabolic subgroup $H \in \Par(\Gamma)$, let 
 \[m_H := \min_{\alpha \in \Par(H)} \{Mod(C_q(\alpha))\} \quad \textrm{and} \quad
 m'_H := \max_{\alpha \in \Par(H)} \{Mod(C_q(\alpha))\} \]
 where $q\in\partial U(H)$.
 Define 
 \[m_\Gamma := \inf_{H\in\Par(\Gamma)} m_H \quad \textrm{and} \quad
  m'_\Gamma := \sup_{H\in\Par(\Gamma)} m'_H.\]
 Since there are finitely many parabolic subgroups up to conjugation, it follows that
 \[0 < m_\Gamma \leq m'_\Gamma < \infty.\]

 \begin{lemma}\label{lem:winding}
 There exists a constant $\sK = \sK(\Gamma)$ such that the following holds.
 Let $H \in \Par(\Gamma)$ be a parabolic subgroup and suppose $\alpha \in \ParC(H)$.
 Then for all $q, q' \in \Disc(\Gamma)$, we have    \[\dist_H(q, q') \asymp_\sK \dist_\alpha(q,q').\]
 \end{lemma}

 \proof
 If $\G$ is a Teichm\"uller geodesic on $\Disc(\Gamma)$ orthogonal to $\partial U(H)$ then
 $\eta_\alpha(q)$ can be chosen to be the same topological leaf for all $q \in \G$.
 Therefore, it suffices to prove the desired result for pairs of points on $\partial U(H)$.
 Fix some $q_0\in \partial U(H)$ and assume, by applying a rotation, that $\alpha$ is horizontal on $(S,q_0)$.
 Then every point on $\partial U(H)$ has the form $q_s = u_s \cdot q_0$ for some $s \in \R$,
 where $u_s$ is the unipotent flow on $\QD(S)$.

 We first estimate the relative twisting under the unipotent flow.
 For brevity, write $\eta_s$ for $\eta_\alpha(q_s)$.
 Observe that $\eta_0$ is vertical on $(S,q_0)$, and has slope $\frac{1}{s}$ on $(S,q_s)$.
 Let $C$ be the cylinder with core curve $\alpha$ on $(S,q_s)$,
 and let $m$ be its modulus.
 Equip the annular cover $S^\alpha$ with the metric obtained by pulling back the half-translation structure
 from $(S,q_s)$, and let $\tilde C$ be the unique Euclidean cylinder on $S^\alpha$ projecting to $C$.
 The number of intersections between $\tilde \eta_0$ and $\tilde \eta_s$ that occur on $\tilde C$
 is equal to $|ms|$ up to a uniform additive error;
 whereas the number of intersections occuring outside of $\tilde C$ is at most two (this follows, for example, using the Gauss--Bonnet Theorem).
 Since $m_\Gamma \leq m \leq m'_\Gamma$, it follows that
  \[m_\Gamma|s| \prec \tw_\alpha(q_0, q_s) \prec m'_\Gamma|s|. \]

 Next, we estimate the cusp winding. Identify $\Disc(q)$ with the upper half-plane model of $\Hy^2$
 so that $\partial U(H)$ coincides with the horizontal line $\Im(z)=1$. 
 Under this identification, the unipotent flow acts as $u_s(z) = z + s$ for all $z \in \Hy^2$.
 In particular, we have $q_s = q_0 + s$ (viewed as points on the complex plane).
 Since $q_s$ is the closest point projection of itself to $\partial U(H)$, we deduce that
 \[\dist_H(q_0, q_s) \asymp |s|.\]
 Combining the above estimates with Proposition \ref{prop:rafi_twist} completes the proof.
 \endproof

 \section{Quasi-isometric embeddings}\label{sec:qi}

 We are now ready to prove the main results.

 \begin{theorem}\label{thm:distance}
  The short marking map    $\smark \from \trCore(\Gamma) \to \Mark(S)$ is a $\Gamma$--equivarient quasi-isometric embedding.
  Furthermore, there exists a constant $c_\Gamma > 0$ such that for all $c \geq c_\Gamma$,
  there exists some $\sK = \sK(S, \Gamma, c)$ such that
  \[\trdist_{\Core(\Gamma)}(x, y) ~\asymp_\sK~ \dist_S(\mu_x, \mu_y)~ + \sum_{\alpha \in \ParC(\Gamma)} \left[\dist_\alpha(\mu_x, \mu_y)\right]_c \]
  for all $x,y \in \trCore(\Gamma)$.
 \end{theorem}

 Since $\smark \from \T_{\geq \epsilon}(S) \to \Mark(S)$ is a $\mcg(S)$--equivariant quasi-isometry, the following is immediate.

 \begin{corollary}\label{cor:thick}
  For $\epsilon > 0$ sufficiently small, the inclusion $\trCore(\Gamma) \hookrightarrow \T_{\geq \epsilon}(S)$ is a quasi-isometric embedding. \qed
 \end{corollary}

 The next result answers a question of Leininger, who originally posed it in the case of
 electrified Teichm\"uller discs arising from lattice Veech groups.
 Our result holds more generally for all finitely generated Veech groups.

 Let us recall the construction of an electrified space.
 Given a length space $\X$ and a collection $\mathcal{U}$ of non-empty subsets of $\X$,
 the \emph{electrification} $\X^{el}$ of $\X$ along $\mathcal{U}$ is defined as follows.
 For each $U \in \mathcal{U}$, introduce a new point $*_U$, called an \emph{electrification point}.
 Then for each $x\in U$, add an interval of length $\frac{1}{2}$ connecting $x$ to $*_U$.
 The metric on $\X^{el}$ is declared to be the induced path metric.
 This procedure forces each subset $U \in \mathcal{U}$ to have diameter at most 1 in $\X^{el}$.
 Note that the inclusion $\X \hookrightarrow \X^{el}$ is 1--Lipschitz.

 Define the \emph{electrified Nielsen core} $\elCore(\Gamma)$ to be
 $\Core(\Gamma)$ electrified along the collection
 of horodiscs $\set{U(H) \st H\in\Par(\Gamma)}$. We shall use $*_H$ to denote the electrification
 point $*_{U(H)}$.
 The systole map can be extended to $\elCore(\Gamma)$ by declaring 
 $\sys(x) := \ParC(H)$ for all $x$ lying in the the open $\frac{1}{2}$--neighbourhood
 of $*_H$ in $\elCore(\Gamma)$. 

 \begin{theorem}\label{thm:electric}
  The systole map $\sys \from \elCore(\Gamma) \to \C(S)$ is a quasi-isometric embedding.
 \end{theorem}

 This result has a natural counterpart for the embedding of the electrified core into
 electrified Teichm\"uller space. For a sufficiently small $\epsilon > 0$,
 define $\elTeich(S)$ to be the electrification of $\T(S)$ along the collection of thin regions
 $V(\alpha) := \set{x \in \T(S) \st \Ext_x(\alpha) \leq \epsilon}$,
 where $\alpha$ runs over $\C(S)$. Write $*_\alpha$ for the associated electrification point.
 The natural inclusion $\iota \from \Core(\Gamma) \to \T(S)$ can be extended to an
 embedding between their respective electrifications as follows.
 By choosing $\epsilon_\Gamma$ and $\epsilon$
 appropriately, we have $U(H) \subseteq V(\alpha) \cap \Core(\Gamma)$
 whenever $H \in \Par(\Gamma)$ and $\alpha \in \ParC(H)$.
 For each $H \in \ParC(\Gamma)$, choose some $\alpha(H) \in \ParC(H)$
 then set $\iota(*_H) = *_{\alpha(H)}$.     We then define $\iota$ on each interval connecting
 $*_H$ to some $x \in U(H)$ in $\elCore(\Gamma)$ by mapping it to the unique interval
 connecting $x$ to $*_{\alpha(H)}$ in $\elTeich(S)$.
 The map $\iota \from \elCore(\Gamma) \to \elTeich(S)$ clearly depends on the choice of cylinder
 curve for each parabolic subgroup, however, it is well-defined up to bounded error.

 The systole map $\sys \from \T(S) \to \C(S)$ can be extended to $\elTeich(S)$ by
 declaring $\sys(x) := \alpha$ for all $x$ in the open $\frac{1}{2}$--neighbourhood
 of $*_\alpha$. Thus, $\sys \circ \iota$ and $\sys$ coarsely agree as maps from $\elCore(\Gamma)$ to $\C(S)$.

 \begin{theorem}[\cite{MM1}]
  The map $\sys \from \elTeich(S) \to \C(S)$ is a quasi-isometry. \halmos
 \end{theorem}

 \begin{corollary}\label{cor:elteich}
  The map $\iota \from \elCore(\Gamma) \to \elTeich(S)$
  is a quasi-isometric embedding. \halmos
 \end{corollary}

 The specific choice of horodiscs is not crucial for the above theorems,
 so long as they are chosen in a $\Gamma$--equivariant manner.

 \subsection{Relevant subsurfaces}

 We require a technical result in order to control the number of terms appearing in sums for our arguments in the following subsections.

 \begin{theorem}[\cite{MM2}]\label{thm:relevant}
  There exists a constant $c_0 = c_0(S) > 0$ for which the following holds.
  Given any threshold $c \geq c_0$, there exists $\sN > 0$ such that
  for any subsurface $Y \subseteq S$ and $\alpha, \beta \in \C(S)$,    the poset
  \[\mathcal{R}^Y_c(\alpha, \beta): = \set{Z \subsetneq Y \st \dist_Z(\alpha, \beta) \geq c} \]
  has at most $\sN \cdot \dist_Y(\alpha, \beta) + \sN$ maximal elements (with respect to subsurface inclusion). \halmos
 \end{theorem}

 This result is a consequence of the Large Links Theorem and the Existence of Hierarchies from from Masur--Minsky \cite{MM2}.
 Their original statement only asserts the existence of a fixed threshold $c_0$
 for which the above holds. We shall sketch a proof of the general statement, allowing for a variable threshold
 $c$, assuming the original statement.

 \proof
 Fix a constant $c \geq c_0$.
 By assumption, there exists a constant $\sN_0 = \sN_0(S) \geq 1$ such that for every essential subsurface $Y \subseteq S$,
 the set $\mathcal{R}^Y_{c_0}(\alpha, \beta)$ has at most $\sN_0 \cdot \dist_Y(\alpha, \beta) + \sN_0$ maximal elements.
 Fix a subsurface $Y$ and let $Z$ be a maximal element of $\mathcal{R}^Y_c(\alpha, \beta)$.
 Then there exists a maximal chain
 \[Y = Y_0 \supsetneq Y_1 \supsetneq \ldots \supsetneq Y_k \supsetneq Z\]
 in $\mathcal{R}^Y_{c_0}(\alpha, \beta)$ where $k \leq \xi(Y)$ and
 $Y_i \notin \mathcal{R}^Y_c(\alpha, \beta)$ for each $0<i\leq k$.
 We shall bound the number of possible chains of this form. Since we consider only maximal chains,
 each $Y_{i+1}$ is a maximal element of $\mathcal{R}^{Y_i}_{c_0}(\alpha, \beta)$.
 Also note that
 \[c_0 \leq \dist_{Y_i}(\alpha, \beta) < c\]
 for each $0 < i \leq k$. Applying the assumption, there are at most 
 $\sN_0\cdot \dist_Y(\alpha, \beta) + \sN_0$ maximal elements in $\mathcal{R}^{Y_0}_{c_0}(\alpha, \beta)$;
 while $\mathcal{R}^{Y_i}_{c_0}(\alpha, \beta)$ has at most $\sN_0\cdot c + \sN_0$
 maximal elements for each $i \geq 1$.
 Therefore, by induction, there are at most
 \[(\sN_0\cdot \dist_Y(\alpha, \beta) + \sN_0) (\sN_0\cdot c + \sN_0)^{\xi(Y)}\]
 possible chains of the desired form.
 Setting $\sN := \sN_0 (\sN_0\cdot c + \sN_0)^{\xi(S)}$ completes the proof.
 \endproof

 \subsection{Proof of Theorem \ref{thm:distance}}
 First, observe that the inclusion $\trCore(\Gamma) \hookrightarrow \T(S)$ is 1--Lipschitz.
 Since the map ${\smark \from \T(S) \to \Mark(S)}$
 is coarsely Lipschitz, it follows that
 \begin{eqnarray}
  \dist_{\Mark(S)}(\mu_x, \mu_y) \quad \prec \quad \trdist_{\Core(\Gamma)}(x, y)\label{eqn:lipschitz}
 \end{eqnarray}
 for all $x,y \in \trCore(\Gamma)$.
 Thus, it remains to prove the reverse coarse inequality.

 We shall establish an upper bound for distances in the truncated core in terms of Teichm\"uller
 distance and cusp winding.

 \begin{lemma}\label{lem:detour}
  For all $c \geq 2$ there exists a constant $\sA_3 > 1$ such that
  \[\trdist_{\Core(\Gamma)}(x,y) \leq \sA_3\cdot\dist_{\Disc(\Gamma)}(x,y) + \sum_{H\in\Par(\Gamma)}  \left[\dist_H(x,y)\right]_c\]
  for all $x,y \in \trCore(\Gamma)$. Moreover, the sum contains finitely many terms.
  \end{lemma}

 \proof
 Let $\pi_\Gamma \from \Disc(\Gamma) \to \trCore(\Gamma)$ be the nearest point projection map.
 Let $\G$ be the Teichm\"uller geodesic connecting $x$ to $y$ in $\Disc(\Gamma)$,
 and $\G'$ be its image under $\pi_\Gamma$. The map $\pi_Y$ replaces each (maximal) subsegment of $\G$
 contained in some horodisc $U(H)$ with a detour running along the associated horocycle $\partial U(H)$
 (with the same endpoints as the given subsegment). Our strategy is to bound the length of each
 detour in terms of the length of the original subsegment.
 This will give an upper bound on the length of $\G'$, and hence $\trdist_{\Core(\Gamma)}(x,y)$.

 For notational convenience, when dealing with any particular parabolic subgroup $H \in \Par(\Gamma)$,
 we shall choose an identification of $\Disc(\Gamma)$ with the upper half-plane $\Hy^2$ so that
 $\partial U(H)$ coincides with the horizontal line $\Im(z) =  1$.
 Under this identification, the map $\pi_\Gamma$ restricted to $U(H)$
 is the vertical projection to $\partial U(H)$.
 Furthermore, we have $\dist_H(x,y) = |\Re(x) - \Re(y)|$, viewing $x$ and $y$ as points
 on the complex plane.

 Consider those $H\in\Par(\Gamma)$ where $\dist_H(x,y) < c$.
 By elementary circle geometry, $\G$ lies
 below the horizontal line $\Im(z) = \sqrt{1 + (\frac{c}{2})^2}$.
 The projection $\pi_\Gamma$ restricted to the region
 lying between this line and $\partial U(H)$ is $\sA_3$--Lipschitz,
 where $\sA_3 = \sA_3(c)$.
 Therefore, the length of the (possibly empty) segment $\G \cap U(H)$
 increases by at most a multiplicative factor of $\sA_3$ under $\pi_\Gamma$.

 Now, consider those $H\in\Par(\Gamma)$ where $\dist_H(x,y) \geq c$. 
 Since $c \geq 2$, the geodesic $\G$ intersects $U(H)$ non-trivially.
 This can only occur for at most finitely many $H \in \Par(\Gamma)$ as the
 pairwise distance between distinct horodiscs is at least 1.
 Observe that $\pi_Y(\G \cap U(H))$ lies inside the horizontal
 line segment connecting $\Re(x) + i$ to $\Re(y) + i$ in $\Hy^2$. Therefore,
 the length of $\G \cap U(H)$ increases by an additive factor of at most $\dist_H(x,y)$
 under $\pi_Y$.

 The above two cases account for all possible detours. The desired result follows
 using the fact that $\G$ has length $\dist_{\Disc(\Gamma)}(x,y)$.
 \endproof

 Next, we choose a suffiently large threshold $c_\Gamma > 2$ to ensure that the following all hold:
 \begin{itemize}
  \item $c_\Gamma \geq c_0$ from Theorem \ref{thm:relevant},
  \item $c_\Gamma \geq c_1$ from Theorem \ref{thm:mm_dist},
  \item $c_\Gamma \geq c_2$ from Theorem \ref{thm:rafi_dist}, where we fix $\epsilon_\Gamma$ as the choice of $\epsilon$, and
  \item $c_\Gamma \geq D_\Gamma$ from Proposition \ref{prop:diameter}; this implies that the only proper subsurfaces
  $Y$ with $\dist_Y(x,y) \geq c_\Gamma$ are parabolic annuli.
 \end{itemize}
 Using Lemma \ref{lem:winding}, we may choose a constant $c' = c'(S, \Gamma, c_\Gamma) \geq 2$ so that
 for any $H \in \Par(\Gamma)$ and ${\alpha \in \ParC(H)}$, we have
 $\dist_\alpha(x,y) \geq c_\Gamma$ whenever $\dist_H(x,y) \geq c'$.
 By Lemma \ref{lem:detour}, we have
 \begin{eqnarray}
  \trdist_{\Core(\Gamma)}(x,y) &\leq& \sA_3\cdot\dist_{\T(S)}(x,y) \quad + \sum_{H\in\Par(\Gamma)}  \left[\dist_H(x,y)\right]_{c'} \label{eqn:sum1}
 \end{eqnarray}
 where $\sA_3 = \sA_3(\Gamma, c')$.
 By Lemma \ref{lem:winding}, there exists some $\sK > 0$ such that
 \begin{eqnarray}
  \dist_H(x,y) &\leq& \sum_{\alpha \in \ParC(H)} \left(\sK \cdot \dist_\alpha(x,y) + \sK\right) \label{eqn:sumH}
 \end{eqnarray}
 for all $H \in \Par(\Gamma)$.
 We shall use (\ref{eqn:sumH}) to replace the sum over parabolic subgroups in (\ref{eqn:sum1}) with a sum
 over parabolic annuli, however, we need to control the number of terms that appear
 due to the additive factor of $\sK$.
 As we are considering only those $H \in \Par(\Gamma)$
 where $\dist_H(x,y) \geq c'$, the corresponding $\alpha$--terms must satisfy
 $\dist_\alpha(x,y) \geq c_\Gamma$.
 Thus, all such (annuli with core curve) $\alpha$ are elements of $\mathcal{R}^S_{c_\Gamma}(\mu_x, \mu_y)$;
 moreover, they are maximal by Proposition \ref{prop:diameter}.
 Therefore, by Theorem \ref{thm:relevant}, there are at most
 $\sN \cdot \dist_S(x,y) + \sN$ relevant parabolic annuli appearing in the sum, where $\sN = \sN(S, c_\Gamma)$.
 Consequently, there exists some $\sB = \sB(S,\Gamma)>0$ such that
 \begin{eqnarray}
  \trdist_{\Core(\Gamma)}(x,y) \quad \prec_{\sB} \quad \dist_{\T(S)}(x,y) \quad + \quad  \left[\dist_S(x,y)\right]_{c_\Gamma} 
  \quad + \sum_{\alpha\in\ParC(\Gamma)}  \left[\dist_\alpha(x,y)\right]_{c_\Gamma}. \label{eqn:Bformula}
 \end{eqnarray}
 Since $\trCore(\Gamma)$ is $\epsilon_\Gamma$--thick, we may apply Theorem \ref{thm:rafi_dist}
 and Proposition \ref{prop:diameter} to obtain
 \begin{eqnarray}
  \dist_{\T(S)}(x,y) &\asymp_{\sA_2}& \sum_{Y \not\cong \Ann} \left[\dist_Y(x,y)\right]_{c_\Gamma}
  + \sum_{\alpha \in \C(S)} \left[\log (\dist_\alpha(x,y))\right]_{c_\Gamma} \\
  & = & \left[\dist_S(x,y)\right]_{c_\Gamma} \quad + \sum_{\alpha\in\ParC(\Gamma)} \left[\log (\dist_\alpha(x,y))\right]_{c_\Gamma} 
 \end{eqnarray}
 for some $\sA_2 = \sA_2(S, \epsilon_\Gamma)$.
 Now, observe that
 \begin{eqnarray}
  [\log t]_c + [t]_c \leq 2[t]_c \label{eqn:log}
 \end{eqnarray}
 for all $t \geq 2$.
 Combining (\ref{eqn:Bformula}) -- (\ref{eqn:log}), we deduce that
 \begin{eqnarray}
  \trdist_{\Core(\Gamma)}(x,y) \quad &\prec_{\sB'}& \quad  \left[\dist_S(x,y)\right]_{c_\Gamma}  \quad + \sum_{\alpha\in\ParC(\Gamma)}  \left[\dist_\alpha(x,y)\right]_{c_\Gamma} \label{eqn:B2formula}\\
  & = & \sum_{Y \subseteq S} \left[\dist_Y(x,y)\right]_{c_\Gamma}
 \end{eqnarray}
 for some $\sB' = \sB'(S,\Gamma)$.
 Finally, combining the above with Theorem \ref{thm:mm_dist} and (\ref{eqn:lipschitz}),
 we may conclude that
 \begin{eqnarray}
  \trdist_{\Core(\Gamma)}(x, y) \quad \asymp_{\sA} \quad \dist_{\Mark(S)}(\mu_x, \mu_y)
 \end{eqnarray}
 for some $\sA = \sA(S,\Gamma)$.
 The desired distance formula, where we allow for any threshold $c \geq c_\Gamma$,
 follows immediately using Theorem \ref{thm:mm_dist}.
 \endproof

 \subsection{Proof of Theorem \ref{thm:electric}}

 The proof proceeds in a similar fashion to the previous subsection.

 \begin{lemma}
  The extended systole map $\sys \from \elCore(\Gamma) \to \C(S)$ is coarsely Lipschitz.
 \end{lemma}

 \proof
 Since $\elCore(\Gamma)$ is a path space, it suffices to prove that any set $V \subset \elCore(\Gamma)$
 of diameter at most $\frac{1}{2}$ has uniformly bounded diameter under $\sys$.
 If $V$ is contained in $\Core(\Gamma)$, then this is immediate from the fact that the usual
 systole map $\sys \from \T(S) \to \C(S)$ is coarsely Lipschitz.

 Now suppose otherwise. Then $V$ non-trivially intersects the $\frac{1}{2}$--neighbourhood
 of some electrification point $*_H$, and is thus contained in the 1--neighbourhood of $*_H$.
 Note that distinct electrification points have disjoint 1--neighbourhoods as the
 pairwise distance between horodiscs is at least 1 in $\elCore(\Gamma)$.
 By taking $\epsilon_\Gamma$ smaller if necessary, the set $\sys(x)$ is contained
 in the simplex $\ParC(H) \subset \C(S)$ for all $x \in U(H)$. By definition, the same is
 true for any $x$ in the (open) {$\frac{1}{2}$--neighbourhood} of $*_H$.
 Any other $x\in \elCore(\Gamma)$ in the 1--neighbourhood of $*_H$
 not accounted for in the previous two cases must lie in the $\frac{1}{2}$--neighbourhood
 of $U(H)$ in $\Core(\Gamma)$.
 The desired result follows using the coarse Lipschitz property of the usual systole map.
 \endproof

 It remains to bound $\eldist_{\Core(\Gamma)}(x,y)$ from above by some linear function of $\dist_S(x,y)$
 for all pairs of points $x,y \in \elCore(\Gamma)$.
 Since $\Core(\Gamma)$ is 1--dense in $\elCore(\Gamma)$, it suffices to prove this for
 $x,y\in\Core(\Gamma)$.

 Choose a constant $c' = c'(S, \Gamma, c_\Gamma) \geq 2$, as in the previous section, so that whenever
 $\dist_H(x,y) \geq c'$ for some $H \in \Par(\Gamma)$,
 we have $\dist_\alpha(x,y) \geq c_\Gamma$ for all ${\alpha \in \ParC(H)}$.  
 Let $\G$ be the Teichm\"uller geodesic in $\Core(\Gamma)$ connecting $x$ to $y$.
 Consider the (finite) set of ${H \in \Par(\Gamma)}$ for which $\dist_H(x,y) \geq c'$.
 We shall order this set $H_1, \ldots, H_n$ according to the order in which 
 the horodiscs $U(H_i)$ appear along $\G$.
 Arguing as in the previous section, we may deduce that 
 \begin{equation}
  n \quad \prec_\sN \quad \dist_S(x,y),
 \end{equation}
 where $\sN = \sN(S, c_\Gamma)$.
 Let $x_i$ and $y_i$ be the endpoints of the
 subinterval $\G \cap \overline{U(H)}$, with $x_i$ chosen to be closer to $x$ along $\G$,
 and set $y_0 = x$ and $x_{n+1} = y$.   
 Applying the the no-backtracking property (Lemma \ref{lem:no_backtrack}) to $\G$ at each of the $x_i$ and $y_i$ for $1 \leq i \leq n$, we deduce that
 \begin{eqnarray}
  \sum_{i=0}^n\dist_S(y_i, x_{i+1})  & \leq & \sum_{i=0}^n\dist_S(y_i, x_{i+1}) \quad + \quad \sum_{i=1}^n \dist_S(x_i, y_i) \\
  & \leq & \dist_S(x,y) + 2Cn \\
  & \prec_\sA &\dist_S(x,y),
 \end{eqnarray}
 where $C = C(S)$ and $\sA = \sA(C, \sN)$.   

 Let $\G_i \subseteq \G$ be the subinterval
 connecting $y_i$ to $x_{i+1}$, for $0 \leq i \leq n$.

 \begin{lemma}
  There exists a constant $\epsilon' = \epsilon'(\Gamma, \epsilon_\Gamma, c') > 0$
  such that each segment $\G_i$ is $\epsilon'$--thick.
 \end{lemma}

 \proof
 Suppose $H \in \Par(\Gamma)$ satisfies $\dist_H(x,y) < c'$.
 Using elementary hyperbolic geometry,
 there exists some $r = r(c')$ such that $\G \cap U(H)$ lies in the $r$--neighbourhood
 of $\partial U(H)$.   Therefore, each $\G_i$ is contained in the $r$--neighbourhood of $\trCore(\Gamma)$ in $\Core(\Gamma)$;
 by cocompactness, all such subintervals are $\epsilon'$--thick for some
 $\epsilon' = \epsilon'(\Gamma, \epsilon_\Gamma, c') > 0$.
 \endproof

 Applying Proposition \ref{prop:thick_geod}, we deduce that $\sys \circ \G_i$ is a parameterised quasigeodesic,
 and so
 \begin{eqnarray}
  |\G_i| & \leq & C' \dist_S(y_i, x_{i+1}) + C'
 \end{eqnarray}
 for some $C' = C'(S,\epsilon')$.

 Next, we construct a modified path $\G'$ in $\elCore(\Gamma)$ by replacing each subsegment
 $\G \cap \overline{U(H)}$ of $\G$ with   the path of length 1
 from $x_i$ to $y_i$ passing through $*_{H_i}$.
 Since the horodiscs are pairwise disjoint, this procedure can be done simultaneously for
 all such horodiscs. By construction, the length of $\G'$ is
 \begin{eqnarray}
  |\G'| &=& \sum_{i=0}^n |\G_i| + n.
 \end{eqnarray}

 Finally, combining the inequalities above, we deduce that
 \begin{eqnarray}
  \eldist_{\Core(\Gamma)}(x,y) & \leq & \sum_{i=0}^n |\G_i| \quad + \quad n\\
  & \leq & C' \left( \sum_{i=0}^n\dist_S(y_i, x_{i+1}) \right)
  \quad + \quad C'(n+1) \quad + \quad n \\
  & \prec_{\sA'} & \dist_S(x,y)
 \end{eqnarray}
 for some $\sA' = \sA'(\sA, C', \sN)$.
 \halmos

\bibliography{mybib}		
\bibliographystyle{amsalpha}

\end{document}